%% file: main.tex
\documentclass[opre,nonblindrev]{informs3}
\DoubleSpacedXI 

\usepackage{endnotes}
\let\footnote=\endnote
\let\enotesize=\normalsize
\def\notesname{Endnotes}%
\def\makeenmark{$^{\theenmark}$}
\def\enoteformat{\rightskip0pt\leftskip0pt\parindent=1.75em
  \leavevmode\llap{\theenmark.\enskip}}

\usepackage{natbib}
 \bibpunct[, ]{(}{)}{,}{a}{}{,}%
 \def\bibfont{\small}%

\usepackage{MnSymbol}
\usepackage[inline]{enumitem}
\usepackage{algpseudocode}
\usepackage{verbatim}
\usepackage{bigstrut}
\usepackage{booktabs,multirow}
\usepackage{lscape}
\usepackage{graphicx}
\usepackage{float}
\usepackage{algorithm}
\usepackage{subfig}
\usepackage{mathtools}
\TheoremsNumberedThrough     
\ECRepeatTheorems

\EquationsNumberedThrough    

\begin{document}


\RUNAUTHOR{Xiangyu Zhang and Peter I. Frazier}

\RUNTITLE{Restless Bandits with Many Arms}

\TITLE{Restless Bandits with Many Arms: \\
Beating the Central Limit Theorem}


\ARTICLEAUTHORS{%
\AUTHOR{Xiangyu Zhang}
\AFF{Department of Operations Research and Information Engineering, Cornell University, Ithaca, NY 14850,  \EMAIL{xz556@cornell.edu}} 
\AUTHOR{Peter I. Frazier}
\AFF{Department of Operations Research and Information Engineering, Cornell University, Ithaca, NY 14850,  \EMAIL{pf98@cornell.edu}}
} 

\ABSTRACT{
We consider finite-horizon restless bandits with multiple pulls per period, which play an important role in recommender systems, active learning, revenue management, and many other areas. 
While an optimal policy can be computed, in principle, using dynamic programming, the computation required scales exponentially in the number of arms $N$. Thus, there is substantial value in understanding the performance of index policies and other policies that can be computed efficiently for large $N$.
We study the growth of the optimality gap, i.e., the loss in expected performance compared to an optimal policy, for such policies in a classical asymptotic regime proposed by Whittle in which $N$ grows while holding constant the fraction of arms that can be pulled per period. 
Intuition from the Central Limit Theorem and the tightest previous theoretical bounds suggest that this optimality gap should grow like $O(\sqrt{N})$. Surprisingly, we show that it is possible to outperform this bound.
We characterize a non-degeneracy condition and a wide class of novel practically-computable policies, called fluid-priority policies, in which the optimality gap is $O(1)$.  These include most widely-used index policies.  When this non-degeneracy condition does not hold, we show that fluid-priority policies nevertheless have an optimality gap that is $O(\sqrt{N})$, significantly generalizing the class of policies for which convergence rates are known. We demonstrate that fluid-priority policies offer state-of-the-art performance on a collection of restless bandit problems in numerical experiments.
}

\maketitle

%


\section{Introduction}


We study a stochastic control problem called the finite-horizon restless bandit. In this problem, 
a decision maker controls $N$ Markov processes (colloquially called ``arms'') with known transition kernels and state-dependent rewards.
The arms produce rewards and evolve independently but are coupled through a constraint (the ``budget'' constraint) on the number that can be activated (colloquially, ``pulled'') in each period. 
Subject to this constraint,
the decision-maker seeks to maximize the expected total reward.

This problem arises in various fields. 
For example, when pursing an active learning approach to classifying images with crowd workers \citep{chen2013optimistic}, each image is an arm, asking a worker to label that image corresponds to pulling that arm, and the arm's state is the resulting Bayesian posterior distribution on the corresponding image's class given past noisy labels.
A limited supply of crowd workers imposes constraints on the number of arms that can be pulled per period.
In dynamic assortment optimization \citep{brown2020index}, a sales manager selects products to display for sale subject to limited display space.
Each product generates revenue at an unknown rate, which can be learned from the revenue it generates when it is displayed.
Each arm is a product, pulling an arm corresponds to displaying that product, and the arm's state is the Bayesian posterior distribution on the product's revenue-generation rate.
Problems in target search by unmanned aerial vehicles
\citep{le2006multi, nino2011sensor},
online advertising \citep{gupta2011thompson, scott2010modern, chakrabarti2009mortal},
network communication \citep{liu2009myopic, al2012multi},
and sensor management \citep{hero2011sensor, nino2011sensor, evans2005networked, nino2011sensor} also fit into our framework.

We study a regime in which the number of arms grows large and the per-period budgets grow proportionally. This regime was first studied in \cite{whittle1980multi}
and has been of longstanding theoretical interest. 
Moreover, it is practically important in many settings.
In examples above, crowdsourced labeling is most challenging when there are many images to label, 
and selecting products for display is most challenging when many products are available.

Despite its importance, this regime presents substantial algorithmic difficulties.
While, in principle, one can compute the optimal policy for restless bandit problems via stochastic dynamic programming, the state of this dynamic program includes the state of each arm and so its dimension grows linearly with $N$. Because of the curse of dimensionality \citep{powell2007approximate}, solving this dynamic program requires computation exponential in $N$.

As a result, there has been substantial interest (e.g., \citealt{whittle1980multi, weber1990index, zayas2019asymptotically, hu2017asymptotically, brown2020index}) in developing approximate policies whose performance is provably close to optimal but require computation that does not grow with $N$. Despite, however, substantial interest and effort focusing on this regime, 
current understanding is limited in several important ways.

First, simulation studies show much better performance for large $N$ in some problems than the best existing theoretical results.
Indeed, the tightest existing upper bound on the optimality gap  
(the difference in performance between the optimal policy and an approximate policy) for such policies is $O(\sqrt{N})$, shown by \cite{brown2020index}
(\citealt{zayas2019asymptotically} provides a policy with a slightly weaker bound of $O(\sqrt{N} \log N)$).
Surprisingly, however, simulation studies by \cite{brown2020index} suggest that the true optimality gap in some problems actually does not grow at all with the number of arms and remains constant at $O(1)$.
The proof techniques used by \cite{brown2020index} and \cite{zayas2019asymptotically}, however, rely heavily on the  Central Limit Theorem (CLT), and do not offer a path toward showing a bound tighter than $O(\sqrt{N})$.

Second, existing theoretical results showing bounds on the optimality gap are restricted to specific policies ($o(N)$, $O(\sqrt{N} \log N)$ and $O(\sqrt{N})$, respectively in \citealt{hu2017asymptotically, zayas2019asymptotically, brown2020index}). 
At the same time, one would expect a very wide class of policies would achieve $o(N)$ and $O(\sqrt{N})$ optimality gaps. 


Our work fills these two gaps: 
we propose a broad class of policies, called \textit{fluid-priority} policies,
which generalize the essential characteristics of policies proposed by \cite{brown2020index} and \cite{hu2017asymptotically}.
Addressing the inconsistency between simulation studies and past theoretical results,
we characterize a sufficient condition, which we call ``non-degeneracy'', under which any fluid-priority policy achieves an $O(1)$ optimality gap, strictly better than all previous results. 
The simulation study consistent with an $O(1)$ optimality gap in \cite{brown2020index} satisfies this non-degeneracy condition.
We also address the current literature's lack of generality 
by providing general easy-to-verify sufficient conditions ensuring $o(N)$ and $O(\sqrt{N})$ optimality gaps. 
All fluid-priority policies satisfy these conditions and thus always achieve an $O(\sqrt{N})$ optimality gap. The policies proposed by \cite{hu2017asymptotically} and \cite{brown2020index} also satisfy the sufficient conditions for an $O(\sqrt{N})$ optimality gap and thus our results generalize those in this previous work.  


To achieve such strong performance, fluid-priority policies generalize well-known index policies by classifying an arm's
state into active, neutral and inactive categories. 
This classification is based on a solution to a linear programming (LP) relaxation that has been important in past analyses of the restless multi-armed bandit problem \citep{whittle1980multi, bertsimas2000restless, hawkins2003langrangian}. 
To be called a fluid-priority policy, 
it should first pull as many arms as possible in active states, up to the budget constraint on the number of arms that can be pulled in this period.
Then, if budget remains, it should should pull arms in neutral states in proportions determined by the solution of the relaxed problem.
Finally, only if budget remains, it should pull arms in inactive states.
There exist many fluid-priority policies because they may prioritize arms in different orders within active, neutral and inactive categories.

Understanding that fluid-priority policies all have good asymptotic performance brings several benefits.
First, it provides a unified understanding of the convergence properties of existing methods, like those proposed in \cite{hu2017asymptotically} and \cite{brown2020index}. 
Second, it can serve as a guideline when developing new policies: it is reasonable to restrict policy development to those within the fluid-priority class.
Third, it creates an opportunity for focused simulation-based search to create policies with excellent empirical performance and provably state-of-the-art asymptotic performance. 
Fluid-priority policies are parameterized by the priority order over active categories, neutral categories, and inactive categories.  (Also, if there are multiple optimal solutions to the relaxed problem, they are additionally parameterized by the choice of solution.)
While still large in problems whose arms have many states, one can perform a focused simulation-based search over this class to find policies with good performance in a specific problem of interest. 
In problems where the number of single-arm states is small enough, it is even possible to search exhaustively over all fluid-priority policies. 
In one numerical experiment, we use this strategy to develop a new fluid-priority policy that significantly outperforms the existing state of the art.
We demonstrate and illustrate these contributions via numerical experiments.
Our first experiment is a Bayesian multi-armed bandit problem with Bernoulli rewards.
We first verify numerically that this problem is non-degenerate.
We then use simulation to calculate expected performance under a fluid priority policy similar to the policies in \cite{hu2017asymptotically} and \cite{brown2020index} and observe that the optimality gap stays constant.
In contrast, we show that the widely used UCB \citep{agrawal1995sample} and Thompson Sampling \citep{agrawal2012analysis} policies have $\Omega(N)$ optimality gaps and significantly underperform by our fluid-priority policy.
Our second experiment is an active learning problem in which one seeks to allocate crowdsourcing effort to accurately classify items, previously studied by \cite{chen2013optimistic}.
We iterate over all possible fluid-priority policies and choose the one with best performance. We find that this fluid priority policy significantly outperforms two previously proposed policies: the Knowledge Gradient \citep{frazier2008knowledge} and Optimistic Knowledge Gradient \citep{chen2013optimistic} policies. Finally, we verify numerically that our non-degeneracy condition holds for dynamic assortment problem studied in \cite{brown2020index}, thus explaining why their simulation study shows an $O(1)$ optimality gap.


Below, we summarize our contribution after first reviewing the literature.

\subsection{Literature Review}
Here we review in more detail the three streams of literature most related to our work. 



\paragraph{Frequentist Bandits:}
The most well-known stream of related work uses frequentist analysis and focuses on problems in which we have uncertainty about an arm's underlying state. 
In such problems, arms are characterized by some underlying but unknown distribution over rewards. This is typically assumed fixed \citep{lai1985asymptotically, auer2002finite}, but can change in some recent analysis \citep{besbes2014stochastic, zhou2020regime}. This literature designs strategies that minimize worst-case expected regret, i.e., the expected difference in total reward compared with a policy that knows arms' underlying characteristics. 

This work is quite different from ours for two reasons. First, it studies a different model using a different performance measure. The model we study assumes that all arms have a fully observed state that evolves stochastically according to Markov processes with known transition kernels and known state-dependent rewards. To apply our model to systems whose arms have unknown reward distributions (which can either be static or vary with a stochastically varying and observable state), one first proposes a Bayesian prior probability distribution over the parameters of these distributions. Then, the Bayesian posterior (which is fully observable) is computed and included as part of the arm's state. Rather than worst-case expected regret, we maximize average case expected reward where initial arm characteristics are drawn at random from the prior. This follows the practice used in partially observable Markov Decision Processes of studying average case reward under a prior belief state. 
The model we study can also be profitably applied to dynamic systems without uncertainty about arms, such as the allocation of airplanes to maintenance bays \citep{cho2015maintenance}.

Second, this difference in model and performance measure creates significant  differences in achievable performance. As we show, policies exist whose average case optimality gap is $O(1)$ in the number of arms $N$. In frequentist bandits, however, the (worst-case) regret grows linearly with $N$ in the simplest case in which arms' characteristics do not change over time \citep{lai1985asymptotically}. 
There is some work that imposes constraints on the relationships between arms, which allows regret to be $o(N)$, such as work on linear bandits \citep{goldenshluger2013linear}, but these models are quite different from the one we consider.

Rather than focusing on $N$, most of this literature focuses on the regime where the horizon $T$ increases to infinity with the number of arms $N$ fixed. 
\cite{lai1985asymptotically} bounds the regret below by a factor proportional to $\log(T)$. Celebrated algorithms such as upper confidence bound (UCB) \citep{auer2002finite} and Thompson Sampling \citep{agrawal2012analysis} are proved to achieve this lower bound asymptotically. 
This stream of work relies on the fact that a long horizon permits a large number of pulls per arm, which distinguish the ``best'' arm from others with high probability.
In our setting, however, where the number of arms is large enough to permit only a small number of pulls per arm and the horizon remains fixed, asymptotic guarantees focusing on large $T$ may not be relevant.
Thus, although there is a large literature demonstrating that variants of UCB \citep{auer2002finite}, Thompson sampling \citep{agrawal2012analysis}, epsilon greedy \citep{sutton1999reinforcement}, and other related algorithms have provably small regret in the large $T$ setting, these results do not imply good performance in the large $N$ setting that we study.
Indeed, in our simulation study, we show that the optimality gap grows linearly with $N$ under UCB and Thompson sampling in the Bayesian finite-horizon multi-armed bandit with Bernoulli reward.

\paragraph{Fixed $N$ and $T$:}
A second and more relevant stream of work \citep{guha2007approximation, guha2008sequential, guha2013approximate, guha2010approximation, farias2011irrevocable, bertsimas2000restless} considers the same model that we consider here and focuses on average-case performance, but considers a regime with a fixed horizon and a fixed number of arms.
This work often solves a linear relaxation of the original problem, constructs a policy based on the solution, and then proves that this policy provides a constant-factor approximation to the optimal one. For example, \cite{farias2011irrevocable} shows that the heuristic they propose achieves an 8-approximation of the optimal policy. Nevertheless, the optimality gap of such policies may scale linearly with $N$, as a constant factor approximation does not preclude this possibility.

\paragraph{Large $N$, fixed $T$:}
The third and most closely related stream considers the same model as the one we consider here and the same asymptotic regime, where 
the number of arms and the budget per period increase proportionally to infinity while holding the horizon fixed. The regime was first studied by \cite{whittle1980multi} in the infinite-horizon discounted reward setting. 
\cite{whittle1980multi} introduced a time-homogeneous Lagrangian relaxation of the budget constraints and proposed the so-called ``Whittle index'' policy when arms are ``indexable'', and conjectured the Whittle index achieves an $o(N)$ optimality gap when this indexability condition holds.
However, \cite{weber1990index} later showed that even under indexability, the optimality gap under the Whittle index policy grows linearly in $N$ for some problems.
Though intuitively promising, the Whittle index policy suffers from the difficulty of verifying the indexibility condition,
the inability to use the policy if indexability does not hold, 
and, in some problems, from weak empirical performance.
Nevertheless, as a pioneering work in restless bandits, the Whittle index inspired a stream of follow-up work, in both the infinite-horizon \citep{bertsimas2000restless, glazebrook2006some, dayanik2008index} and finite-horizon cases. As it is the focus of our work, we now discuss the finite-horizon case in detail.


Following Whittle's earlier work, later literature (e.g., \citealt{hu2017asymptotically, zayas2019asymptotically, brown2020index}) studies the finite-horizon restless bandit using Lagrangian relaxations. Unlike Whittle's work, these use time-dependent Lagrange multipliers because of the non-stationary nature of finite-horizon problems.
This technique yields promising performance guarantees and empirical results without the need for an indexability condition.
\cite{hu2017asymptotically} studies the binary-action bandit problem and proposes an index policy achieving an $o(N)$ optimality gap. \cite{zayas2019asymptotically} studies the multi-action bandit problem and proposes a policy achieving an $O(\sqrt{N} \log N)$ optimality gap.
\cite{brown2020index} studies the same setting as \cite{hu2017asymptotically} and proposes policies with an $O(\sqrt{N})$ optimality gap. However, simulation experiments \citep{brown2020index} suggest, surprisingly, that the optimality gap might not grow with $N$. Our work proposes a novel policy class, the class of {\it fluid-priority} policies, and shows that policies proposed by \cite{hu2017asymptotically, brown2020index} are special cases in this class. Furthermore, we show that any policy in this class achieves an $O(\sqrt{N})$ optimality gap in all circumstances and achieves an $O(1)$ optimality gap when a non-degeneracy condition holds. Specially, for the setting discussed above in which simulation experiments from \cite{brown2020index} suggest the optimality gap is $O(1)$, we show the non-degeneracy condition holds.

\subsection{Summary of Contributions and Outline}
There are three main contributions in our work.

\textit{Main Contribution}: Our main contribution is to propose a novel and general class of policies, fluid-priority policies, and show that they have strong theoretical performance guarantees. We show theoretically that any fluid-priority policy achieves an $O(\sqrt{N})$ optimality gap in all circumstances and achieves an $O(1)$ optimality gap under a non-degeneracy condition.

\textit{Secondary Contributions}: Building on our main result, we provide three secondary contributions.
\begin{enumerate}
\item We establish $o(N)$ and $O(\sqrt{N})$ optimality gaps for classes of policies broader than fluid priority policies. The sufficient conditions used by our proof are easy to verify and general enough to apply to policies proposed in \cite{hu2017asymptotically, brown2020index}. 

\item We propose an algorithm for verifying whether non-degeneracy condition holds. 
If so, this algorithm searches over optimal occupation measures to find one that is non-degenerate.

\item We demonstrate the value of fluid policies through numerical experiments and additional theory.
\begin{itemize}
\item By searching numerically over fluid priority policies, we identify a novel fluid priority policy that outperforms a previously proposed state-of-the-art policy designed specifically for crowdsourced labeling. 
\item We show that the dynamic assortment problem previously studied by \cite{brown2020index} satisfies the non-degeneracy condition, and thus our theoretical results explain the hitherto poorly understood performance of Lagrangian index policies in this setting. 
\item We show theoretically that the widely-used UCB and Thompson sampling algorithms for finite-horizon Bernoulli bandits have strictly worse asymptotic expected performance (in the worst case over problem instances) than fluid priority policies for problems with many arms. We then demonstrate numerically that fluid priority policies have substantially better empirical performance in a collection of such problems.
\end{itemize}
\end{enumerate}

\textit{Organization of This Paper}:
The rest of the paper is organized as follows.
\S\ref{sys_model} defines the restless bandit problem as a MDP.
\S\ref{approx} introduces notation and provides background on an existing linear programming relaxation used in our later novel theoretical results. This relaxation provides an upper bound on the problem's optimal performance.
Based on the upper bound, \S\ref{oN} describes sufficient conditions to achieve an $O(N)$ optimality gap and \S\ref{sec-4} describes sufficient conditions to achieve an $O(\sqrt{N})$ optimality gap. 
\S\ref{sect-6} proposes the class of fluid-priority policies and proves that they achieve an $O(\sqrt{N})$ optimality gap. \S\ref{se-7} proves that fluid-priority policies achieve an $O(1)$ optimality gap when the non-degeneracy condition is met. \S\ref{numeric} provides numerical studies and \S\ref{conclusion} concludes our work.

\section{System Model} \label{sys_model}
This section formulates our decision-making problem as a Markov Decision Process (MDP). 

\textbf{Model:}
There are $N$ arms, each of which shares the same finite state space $S$. We use $s_{i, t}$ to indicate the state of arm $i$ at time $t$. At each period $t$ for each arm $i$, the decision-maker chooses whether to pull the arm ($a_{i,t} = 1$) or leave it idle ($a_{i,t} = 0$). We define $A = \{0, 1 \}$ to be the space of available actions in which $a_{i,t}$ takes values. 
These actions must respect a so-called ``budget constraint'' in which the number of arms pulled at period $t$ is $B_t = \lfloor \alpha_t N \rfloor$, where $0 \leq \alpha_t \leq 1$ is a pre-specified budget ratio.

Based on the action applied, each arm's state transitions stochastically to time $t + 1$ according to a known transition kernel
    $P_t = \{ p_t(s, a, s') \}_{s, s' \in S, a \in A}$ 
    where
    $p_t(s, a, s') = \mathbb{P}(s_{t + 1, i} = s' | s_{t, i} = s, a_{t, i} = a)$.
All arms share the same transition kernel, and any arm's transition is conditionally independent from others given its own state and action. (Arm-specific transition kernels can be modeled by defining static arm ``types'' and extending the state space to specify the arm's type.) 
At time period $t$, each state-action pair is associated with a reward, given by a known reward function $r_t: S \times A \rightarrow \mathbb{R}$.  The decision-maker aims to maximize the total reward collected from all $N$ arms over a finite horizon subject to the budget constraint.

To complete the formal definition of our problem involving $N$ arms, we introduce some additional notation. We use $\mathbb{S} = S^N$ to denote the $N$-fold Cartesian product of the state space $S$ and define $\mathbb{A} = A^N$ similarly.
All $N$ arms together form an MDP with state space $\mathbb{S}$ and action space $\mathbb{A}$. We call this the ``joint MDP'' to distinguish it from MDPs that we reference later involving a single arm. 
The state in this joint MDP at time $t$ is 
$\mathbf{s}_t = (s_{t, 1}, s_{t, 2}, ..., s_{t, N}) \in \mathbb{S}$, which indicates that arm $i$ has state $s_{t, i}$.
The action is $\mathbf{a}_t = (a_{t, 1}, a_{t, 2}, ..., a_{t, N}) \in \mathbb{A}$, which indicates that action $a_{t, i}$ is applied to arm $i$.

The reward function of the joint MDP, $R_t: \mathbb{S} \times \mathbb{A} \rightarrow \mathbb{R}$, is the sum of the single-arm rewards defined above,
\begin{equation*}
R_t(\mathbf{s}_t, \mathbf{a}_t) = \sum_{i=1}^N r_t(s_{t, i}, a_{t, i}).
\end{equation*}

For element $\mathbf{a} = (a_{1}, a_{2}, ..., a_{N})$ in $\mathbb{A}$, we use $|\mathbf{a}| = \sum_{i=1}^N a_i$ to indicate the $L^1$-norm of $\mathbf{a}$, i.e, the number of pulled arms. We write our budget constraint at time $t$ as $|\mathbf{a}_t| = B_t$.


The transition kernel for the joint MDP is the product of each arm's transition kernel,
\begin{align*}
    \mathbb{P}[\mathbf{s}_{t + 1} | \mathbf{s}_t, \mathbf{a}_t] = \prod_{i = 1}^N p_t(s_{t, i}, a_{t, i}, s_{t + 1, i}).
\end{align*}

We assume all arms start from the same initial state $s^*$.
Our analysis can be easily generalized to the case where arms start from different states.

A policy $\pi$ is a function that maps the current state $\mathbf{s}_t \in \mathbb{S}$ and time $t$ to an action $\mathbf{a}_t \in \mathbb{A}$. The objective of the policy is to maximize the expected total reward, subject to the budget constraint specified above.

This objective can be written as,
\begin{equation}\label{OP}
\begin{split}
    &\max_{\pi} \ \mathbb{E}_{\pi} \sum_{t=1}^T R_t(\mathbf{s}_t, \mathbf{a}_t) \\
    &\text{subject to:} \ |\mathbf{a}_t| = \lfloor \alpha_t N \rfloor, \ \forall t \in [T],
\end{split}
\end{equation}
where $\mathbb{E}_{\pi}$ indicates the expectation taken under policy $\pi$.

We define the value function of a policy $\pi$ as $V_N(\pi) = \mathbb{E}_{\pi} \sum_{t=1}^T R_t(\mathbf{s}_t, \mathbf{a}_t)$. We measure a policy's performance by comparing its value with that of an optimal policy solving (\ref{OP}). Let $V_N^* = \sup_{\pi} V_N(\pi)$
be the value of an optimal policy. Then the optimality gap of the policy $\pi$ is defined as
\begin{align*}
    V_N^* -  V_N(\pi).
\end{align*}
Maximizing the value function across policies is equivalent to minimizing the optimality gap. We are interested in finding policies with small optimality gaps when $N$ is large.

\textbf{Applications:}
The above model has many applications. 
In the most direct application, each arm corresponds to a physical process that evolves stochastically and independently of the other physical processes according to a known transition kernel. 
Examples include network communication \citep{liu2009myopic, al2012multi}  and machine maintenance \citep{glazebrook2006some, abbou2019group, cho2015maintenance}. 
For example, in maintenance of military aircraft 
with low radar visibility (so-called ``stealth'' aircraft) \citep{cho2015maintenance}, each aircraft is treated as an arm. Radar visibility (the state of the arm) increases stochastically according to a known transition kernel each time the aircraft flies as small particles in the air damage the aircraft's paint and underlying metal surface. This damage can be repaired (the arm can be pulled) by pausing an aircraft's flights and performing maintenance. Our objective is allocate limited maintenance resources to maximize an objective combining flights flown and number of aircraft with low radar visibility.


In addition, there are many applications in which {\it information} evolves over time. In such settings, we often have several independent unknown quantities, each arm corresponds to one of these quantities, and an arm's state represents the information that we have about this quantity. Examples include autonomous target tracking \citep{le2006multi, hero2011sensor},  where each target is treated as an arm, and its state is whether it is tracked by a sensor and some physical feature affecting the motion of the target. Based on its state, the target moves to a new location, and our objective is to track as many targets for as long as possible.

In perhaps the most famous restless bandit, each arm corresponds to a slot machine. Each slot machine generates payoffs according to a distribution from a parametric family (e.g., Bernoulli). The parameter governing an arm's rewards (for Bernoulli arms, the payoff probability) is drawn at random from a Bayesian prior distribution and is unobserved. The state of the arm is the Bayesian posterior distribution on its parameter, given all observed payoffs from the arm. When we pull an arm, we earn a reward (whose distribution is given by marginalizing over the posterior on the arm's uncertain parameter) and the new state is determined by Bayes' rule and the observed reward. If an arm's underlying parameter changes over time, then this causes the posterior to change even if the arm is not pulled, making the problem restless. A common point of confusion arises from the fact that this problem has a similar motivation to the more widely studied non-Bayesian stochastic bandit \citep{lai1985asymptotically, auer2002finite}, but uses a different formulation. For a tutorial on Bayesian formulations of multi-armed bandits, some of which are restless, see \cite{mahajan2008multi}.

\section{Background: Preliminary Results and Notation}
\label{approx}

In this section, we define a linear programming relaxation that bounds $\hat{V}_N^*$ for $V_N^*$. 
Although this bound is standard in the literature and is not part of our contribution, we include it to provide a self-contained presentation and to establish notation used later.


\textbf{Linear Programming Relaxation:}
Similar to \cite{wu2015algorithms, farias2011irrevocable, guha2008sequential}, we introduce this relaxation of Problem (\ref{OP}):
\begin{equation}\label{KP}
\begin{split}
    & \hat{V}_N^* := \max_{\pi} \ \mathbb{E}_{\pi} \sum_{t=1}^T R_t(\mathbf{s}_t, \mathbf{a}_t) \\
    &\text{subject to } \mathbb{E}_{\pi}|\mathbf{a}_t| = \alpha_t N, \ \forall t \in [T].
\end{split}
\end{equation}
This relaxes problem (\ref{OP})'s almost sure cardinality constraints (on both the initial occupation measure and the number of pulls) to constraints on the expected cardinality.
As we will see soon, solving relaxation (\ref{KP}) is equivalent to solving a linear program whose number of decision variables does not depend on $N$ (see Lemma \ref{fenchel} and the linear program (\ref{LP})).
For simplicity of presentation, we assume that
$\alpha_t$ are rational and we restrict attention and limits taken below over $N$ causing $\alpha_t N$ to be integral for all $t \in [T]$.
Our results essentially generalize to irrational $\alpha_t$ and non-integral $\alpha_t N$ as discussed briefly in Appendix \ref{al:2-1}.

The value of this relaxed problem, $\hat{V}_N^*$, is an upper bound on $V_N^*$. We use this upper bound extensively later to bound the optimality gap of the policies we study. Moreover, the policies we study in \S\ref{sect-6} heavily leverage this relaxation in their definition.  They benefit from the fact that the relaxation yields a low-dimensional problem whose number of decision variables and constraints do not scale with $N$. This allows the relaxation's solution to be computed and used to define practical policies, even when $N$ is large. 




The following lemma formally states this bound and also observes (via Fenchel's duality theorem, and the separability of a dualized version of Problem (\ref{KP})) that $\hat{V}_N^*$ is determined by the solution to a single-armed problem $\hat{V}^*_1$.
Its proof can be found in the appendix.
\begin{lemma}\label{fenchel}
$V_N^* \leq \hat{V}_N^* = N \hat{V}_1^*.$
\end{lemma}

The quantity $\hat{V}_1^*$ is given by, 
\begin{equation}\label{P}
\begin{split}
    &\max_{\pi} \ \mathbb{E}_{\pi} \sum_{t=1}^T r_t(s_t, a_t) \\
    &\text{subject to } \mathbb{E}_{\pi}|{a}_t| = \alpha_t, \ \forall t \in [T].
\end{split}
\end{equation}

Later analysis and computation is supported by the following equivalent version of Problem (\ref{P}).
Defining the occupation measure, $x_t(s, a) := \mathbb{P}[s_t = s, a_t = a]$, 
Problem (\ref{P}) is equivalent to
\begin{equation}\label{LP}
\begin{split}
    & \max \  \sum_{s \in S, a \in A}\sum_{t=1}^T r_t(s, a) x_t(s, a) \\
    &\text{subject to }\\
    &\quad \quad \quad 
    \sum_{a \in A}x_t(s, a) = \sum_{a \in A} \sum_{s' \in S} x_{t-1}(s', a)p_{t-1}(s', a, s), \text{ $\forall s \in S,  2 \leq t \leq T$, } \\
    &\quad \quad \quad
    \sum_{s \in S}x_t(s, 1) = \alpha_t, \ \text{$t \in [T]$, } \\
    &\quad \quad \quad 
    \sum_{a \in A}x_1(s^*, a) = 1,\\
    &\quad \quad \quad 
    \sum_{a \in A}\sum_{s \in S}x_1(s, a) = 1, \text{ $\forall s \in S$,}  \\
    &\quad \quad \quad
    x_t(s, a) \geq 0, \text{ $\forall s \in S, a \in A,  t \in [T]$.} 
\end{split}
\end{equation}
The first constraint of Problem (\ref{LP}) ensures that flows are balanced; the second ensures that the budget constraint is met; and the third follows from the initial occupation measure. 
We let $x_t(s,a)$ denote the entries in an optimal occupation measure, i.e., one that solves Problem (\ref{LP}). Then, we can compute, 
\begin{align}\label{eq-5}
    \hat{V}^*_1 = \sum_{s \in S, a \in A}\sum_{t=1}^T r_t(s, a) x_t(s, a).
\end{align}

The class of policies we analyze depend on solving Problem (\ref{LP}) computationally using a linear programming solver. As noted above, this is possible, even when $N$ is large, because the dimensionality of Problem (\ref{LP}) does not depend on the number of arms $N$.

\textbf{Additional Notation:}
Here we introduce some additional notation used in the following sections.
Given the optimal occupation measure, we use
$z_t(s) := \sum_{a \in A}x_t(s, a)$ to denote the probability that an arm is in state $s$ at time $t$ under this measure. We use $z_t$ and $x_t$ to refer to the corresponding vector (or matrix), i.e., $z_t := (z_t(s), s \in S)$ or $x_t := (x_t(s, a) : s \in S, a \in A)$.

In the joint MDP with $N$ arms, we let 
$X^N_t(s, a)$ be the number of arms in state $s$ for which we take action $a$ at time $t$.
We let $Z^N_t(s)$ be the number of arms in state $s$ at time $t$.
We use $Z_t^N, X_t^N$ to refer to the vectors $(Z_t^N(s) : s \in S)$ and matrix $(X_t^N(s, a) : s \in S, a \in A)$.
Using this notation, a policy $\pi$ of the joint MDP is a map from $Z^N_t$ to $X^N_t$.

\S\ref{sec-4} will study deviations between the realization of $(Z_t^N, X_t^N)$ and $(N z_t, N x_t)$, and how these deviations impact the joint MDP's reward. To support this analysis, we define {\it{diffusion statistics}} $\tilde{Z}^N_t$ and $\tilde{X}^N_t$ as
\begin{equation*}
    \tilde{Z}^N_t = \frac{Z_t^N - N z_t}{\sqrt{N}},\quad
    \tilde{X}^N_t = \frac{X_t^N - N x_t}{\sqrt{N}}.
\end{equation*}
Using this notation, a policy $\pi$ of the joint MDP naturally induces a class of maps $\tilde{\pi}_{t, N}$ indexed by $t$ and $N$, from diffusion $\tilde{Z}^N_t$ to diffusion $\tilde{X}^N_t$, such that
\begin{align}\label{eq-6}
    \pi(t, Z_t^N) = X_t^N \Longleftrightarrow \tilde{\pi}_{t, N}(\tilde{Z}_t^N) = \tilde{X}_t^N.
\end{align}

\section{Sufficient Conditions for Achieving an $o(N)$ Optimality Gap}\label{oN}
This section establishes the first of our contributions: general sufficient conditions for an $o(N)$ optimality gap.
This result allows us to directly verify that the policy in 
\cite{zayas2019asymptotically} has an $o(N)$ optimality gap.
We build on the results here in the next section, where we give stronger conditions sufficient for an $O(\sqrt{N})$ gap and apply it to the policies in \cite{hu2017asymptotically} and \cite{brown2020index}. This is in preparation for our main contribution in \S\ref{sect-6}, a class of policies with an $O(1)$ gap.

The main idea in this section is, essentially, that as long as the number of arms we pull in each state, $X^N_t$, is approximately proportional to the optimal occupation measure $x_t$ (a property we formalize and give the name ``fluid consistency''), the number of arms in the next period $Z_{t + 1}^N$ in each state will be approximately proportional to $z_{t + 1}$. This will cause the reward of the joint MDP to scale proportionally with $\hat{V}^*_1$. While random fluctuations cause proportionality to hold only approximately, their resulting loss in reward is $o(N)$.

We begin by formally defining the notion of fluid consistency.
\begin{definition}
Under a policy $\pi$, 
if $\pi(t, Z_t^N)/N \to x_t$ for all $t\in[T]$ and sequences $(Z_t^N : N)$ satisfying 
$Z_t^N/N \rightarrow z_t$,
then we say the policy $\pi$ is {\it fluid consistent}.
\end{definition}

Based on this definition, we have the following lemma, whose proof can be found in the appendix.
\begin{lemma}\label{fluid_consistent}
If a policy $\pi$ is fluid consistent, then
\begin{align*}
    \frac{Z_t^N}{N} \rightarrow z_t, \frac{X_t^N}{N} \rightarrow x_t, 
\end{align*}
almost surely for any $t \in [T]$ as $N\to\infty$.
\end{lemma}

Using Lemma \ref{fluid_consistent}, we now show the main result of this section: that fluid consistency implies the optimality gap is $o(N)$.
\begin{theorem}\label{fluid_optimal}
If a policy $\pi$ is fluid consistent, then $V_N^* - V_N(\pi) = o(N)$.
\end{theorem}

\proof{Proof of Theorem \ref{fluid_optimal}}
Because the policy $\pi$ is fluid consistent, Lemma \ref{fluid_consistent} shows
\begin{align*}
    \frac{Z_t^N}{N} \rightarrow z_t, \frac{X_t^N}{N} \rightarrow x_t.
\end{align*}

The total reward of the joint MDP, divided by $N$, is
\begin{align*}
    \frac{1}{N}\mathbb{E}_{\pi} \sum_{t = 1}^T R_t(\mathbf{s}_t,\mathbf{a}_t)
    &= \frac{1}{N} \mathbb{E}_{\pi} \sum_{t = 1}^T \sum_{s \in S, a \in A} r_t(s, a) X_t^N(s, a)\\
    &= \mathbb{E}_{\pi} \sum_{t = 1}^T \sum_{s \in S, a \in A} r_t(s, a) \frac{X_t^N(s, a)}{N}\\
    &\rightarrow \mathbb{E}_{\pi} \sum_{t = 1}^T \sum_{s \in S, a \in A} r_t(s, a) x_t(s, a)
\end{align*}
as $N\to\infty$, where we leverage the dominated convergence theorem, the fact that rewards are bounded, and $0 \le X_t^N(s, a) \le N$.
Thus, we have shown that $V_N^* - V_N(\pi) = o(N)$.
\hfill $\Box$

One can show that the the policies in \cite{hu2017asymptotically, zayas2019asymptotically,brown2020index} are all fluid consistent and thus have $o(N)$ optimality gaps.
We show this for \cite{zayas2019asymptotically} in Appendix \ref{other-policies}.
Below, we show that \cite{hu2017asymptotically, brown2020index} meet a stronger condition and thus have $O(\sqrt{N})$ optimality gaps.

\section{Sufficient Conditions for Achieving an $O(\sqrt{N})$ Optimality Gap}\label{sec-4}
This section establishes our second contribution: a substantially more general result than in the literature showing sufficient conditions for an $O(\sqrt{N})$ optimality gap.
Using this result, we directly verify that policies in \cite{hu2017asymptotically} and \cite{brown2020index} have $O(\sqrt{N})$ optimality gaps.
This section also provides stepping stones towards our main contribution, described in 
\S\ref{sect-6}.

The main idea in this section is that,
as long as the diffusion statistic $\tilde{X}_t^N$ is bounded by $O(1)$, then $\tilde{Z}_{t + 1}^N$ will also be bounded by $O(1)$. Thus, the deviation between the reward of the joint MDP and the relaxation's upper bound $\hat{V}_N^*$ will be bounded by $\sqrt{N}\cdot O(1) = O(\sqrt{N})$.

Recall Equation (\ref{eq-6}), that a policy $\pi$ naturally induces a class of maps $\tilde{\pi}_{t, N}$. Using this idea, we say a policy $\pi$ is ``diffusion regular'' if all induced maps $\tilde{\pi}_{t, N}$ keep the diffusion $\tilde{X}_t^N$ bounded by $O(1)$. We define this formally here.
\begin{definition}\label{def-diff}
A policy $\pi$ is called {\it diffusion regular} if its induced maps $\tilde{\pi}_{t, N}$ satisfy the following conditions, where $|\cdot|$ indicates the $L^1$-norm in Euclidean space.
\begin{enumerate}
  \item There exists $C_1 > 0$ s.t.
    $|\tilde{\pi}_{t, N}(\theta_1) - \tilde{\pi}_{t, N}(\theta_2)| \leq C_1 |\theta_1 - \theta_2|$
  for all $t$, $N$, $\theta_1$ and $\theta_2$.
  \item There exists $C_2 > 0$ s.t.
    $|\tilde{\pi}_{t, N}(0)| \leq C_2$
  for all $t$ and $N$.
    \item There exists a map $\tilde{\pi}_{t, \infty}$ s.t. $\tilde{\pi}_{t, N}(\theta) \rightarrow \tilde{\pi}_{t, \infty}(\theta)$ as $N \rightarrow \infty$ 
    for all $\theta$.
\end{enumerate}
\end{definition}

We briefly note the following fact, useful when proving subsequent results. Its proof is found in the appendix.
\begin{lemma}\label{fluid-imply-diffusion}
If a policy is diffusion regular then it is also fluid consistent.
\end{lemma}

We now show that
if a policy $\pi$ is diffusion regular, the diffusion statistics $\tilde{X}_t^{\infty}$ and $\tilde{Z}_t^{\infty}$ converge in distribution (Lemma \ref{diffusion_consistent}) and their second moments are uniformly bounded (Lemma \ref{moment-bound}). Proofs of Lemma \ref{diffusion_consistent} and Lemma \ref{moment-bound} are in the Appendix.
\begin{lemma} \label{diffusion_consistent}
If a policy $\pi$ is diffusion regular, then for any $t \in [T]$, there exists sub-Gaussian random vectors $(\tilde{Z}_t^{\infty},\tilde{X}_t^{\infty})$ such that
$(\tilde{Z}_t^N, \tilde{X}_t^{N}) \rightarrow (\tilde{Z}_t^{\infty},\tilde{X}_t^{\infty})$
in distribution as $N \rightarrow \infty$. 
\end{lemma}

\begin{lemma} \label{moment-bound}
If a policy $\pi$ is diffusion regular, then there exists a constant $C$ such that
$\mathbb{E}_{\pi}[||\tilde{Z}_t^N||_2^2] \leq C$ and $\mathbb{E}_{\pi}[||\tilde{X}_t^N||_2^2] \leq C$
for all $t \in [T]$ and $N$, where $||\cdot||_2$ indicates the $L^2$ norm.
\end{lemma}

Based on Lemma \ref{diffusion_consistent} and \ref{moment-bound}, we can prove the following theorem.
\begin{theorem} \label{diffusion_optimal}
If a policy $\pi$ is diffusion regular, then $V_N^* - V_N(\pi) = O(\sqrt{N})$.
\end{theorem}

\proof{Proof of Theorem \ref{diffusion_optimal}}
Since the policy $\pi$ is diffusion regular, there exists sub-Gaussian random vectors $\tilde{Z}_t^{\infty}, \tilde{X}_t^{\infty}$, such that 
$\tilde{Z}_t^N \rightarrow \tilde{Z}_t^{\infty}$ and $\tilde{X}_t^{N} \rightarrow \tilde{X}_t^{\infty}$ in distribution as $N\to \infty$ by Lemma \ref{diffusion_consistent}. 

Also, the optimality gap is bounded above by
\begin{align*}
    V_N^* - V_N(\pi)
    &\le N\hat{V}_1^* - V_N(\pi)\\
    &= N\sum_{t = 1}^T \sum_{s \in S, a \in A}r_t(s, a)x_t(s, a) - \mathbb{E}_{\pi} \sum_{t = 1}^T R_t(\mathbf{s}_t,\mathbf{a}_t) \\
    &= -\sqrt{N} \mathbb{E}_{\pi} \sum_{t = 1}^T \sum_{s \in S, a \in A} r_t(s, a) \tilde{X}_t^N(s, a).
\end{align*}

Divide both sides of this bound by $\sqrt{N}$ and take $N \rightarrow \infty$. Then, since $\tilde{X}_t^N$ and $\hat{Y}_t^N$ are uniformly integrable (Lemma \ref{moment-bound}),
\begin{equation*}
    \lim \sup_N \frac{1}{\sqrt{N}} (V_N^* - V_N(\pi)) \leq \lim \sup_N -\mathbb{E}_{\pi} \sum_{t = 1}^T \sum_{s \in S, a \in A} r_t(s, a) \tilde{X}_t^N(s, a)
    = -\mathbb{E}_{\pi} \sum_{t = 1}^T \sum_{s \in S, a \in A} r_t(s, a) \tilde{X}_t^{\infty}(s, a).
\end{equation*}

To summarize, we have shown $V_N^* - V_N(\pi) = O(\sqrt{N})$.
\hfill $\Box$

We verify in the Appendix \ref{other-policies} that the 
policies proposed by \cite{hu2017asymptotically} and \cite{brown2020index} are diffusion regular
and thus (by Theorem \ref{diffusion_optimal}) have $O(\sqrt{N})$ optimality gaps.
Thus, Theorem \ref{diffusion_optimal} generalizes the performance guarantees shown in that previous work.

\section{Fluid-priority policies}\label{sect-6}
This section defines fluid-priority policies and show that they are always diffusion regular and thus achieve an $O(\sqrt{N})$ optimality gap. Later, in Section~7, we show that they achieve an $O(1)$ optimality gap if an additional condition is satisfied.

Roughly speaking, a fluid-priority policy is defined by first fetching an optimal solution of the LP relaxation, then classifying states into three disjoint categories based on the solution: fluid-active, fluid-neutral and fluid-inactive. A fluid priority policy is one that pulls arms while respecting a prioritization derived from these categories: arms in fluid-active states are prioritized for pulling over those in fluid-neutral states; and arms in fluid-neutral states are prioritized in turn over arms in fluid-inactive states.
Additionally, when pulling arms in fluid-neutral states, a fluid-priority policy must do so according to proportions derived from LP relaxation.


Mathematically speaking, a fluid-priority policy is parameterized by an occupation measure $\{x_t(s, a)\}_{t, s, a}$ solving Problem ($\ref{LP}$) and a sequence of ``priority-score'' functions $\{ \mathcal{P}_t(\cdot) \}_t$ assigning each state a real number. 
Based on the occupation measure $\{x_t(s, a)\}_{t, s, a}$, a fluid-priority policy classifies states into these three disjoint categories:
\begin{align}\label{category}
&\text{The {\it fluid-active} category:}\  C^{+}_t := \{ s \in S | x_t(s, 1) > 0, x_t(s, 0) = 0\},\nonumber\\
&\text{The {\it fluid-neutral} category:}\ C^{0}_t := \{ s \in S | x_t(s, 1) > 0, x_t(s, 0) > 0\},\\
&\text{The {\it fluid-inactive} category:}\ C^{-}_t := \{ s \in S | x_t(s, 1) = 0, x_t(s, 0) = 0\}\nonumber.
\end{align}
We refer to an arm whose state is in the fluid-active category as a {\it fluid-active arm}. We define the terminology {\it fluid-neutral arm} and {\it fluid-inactive arm} similarly.
With these definitions in place, the fluid-priority policy corresponding to a given occupation measure and priority-score function is defined by Algorithm~1.





Algorithm \ref{al:1} allocates its budget by first pulling as many fluid-active arms as possible, subject to the budget constraint (Lines 5-7).  If budget remains, then it pulls as many fluid-neutral arms as possible, again subject to the constraint on the remaining budget (Lines 9-17). 

When there is enough budget to pull some fluid-neutral arms, but not all of them, the budget is allocated carefully across them to ensure fluid-consistency.
This is closely related to ``tie-breaking'' as discussed in Algorithm 2 of \cite{hu2017asymptotically}.
In particular, lines 9-13 ensure that the number of arms pulled in each fluid-neutral state is at least equal to $\lfloor N x_t(s,1) \rfloor$, the number of arms from this state pulled in the fluid relaxation, as long as the budget constraint $B_t$ and number of available arms $Z_t(s)$ allows.
If budget remains after this is achieved, additional fluid-neutral arms are pulled. 

Finally, if budget remains after all fluid-neutral arms are pulled, additional fluid-inactive arms are pulled until the budget is exhausted. 
Within each category (fluid-active, fluid-neutral, fluid-inactive), states are prioritized based on the priority score.

\begin{algorithm}
	\caption{Fluid-priority policy}\label{al:1}
	\hspace*{\algorithmicindent} \textbf{Input:} optimal occupation measure $(x_t(s, a))_{t \in [T], s \in S, a \in A}$ found by solving the linear program \eqref{LP} and priority-score functions $\{ \mathcal{P}_t \}_{t \in [T]}$.
	\begin{algorithmic}[1]
	\For {$t = 1, 2, ..., T$}
	    \State Classify states into fluid-active ($C^+_t$), fluid-neutral ($C^0_t$) and fluid-inactive ($C^-_t$) categories based on the occupation measure, according to equation (\ref{category}).
		\State Observe there are $Z_t(s)$ arms in state $s$ and remaining budget $B_t = \lfloor \alpha_t N \rfloor$.
		\State Each of the for loops below iterates over states in decreasing order of $\mathcal{P}_t(s)$
		\For {state s in $C^{+}_t$}\label{1l:9}
    		\State Plan to pull $X_t^N(s, 1) \leftarrow \min\{B_t, Z_t(s)\}$ arms out of the $Z_t(s)$ arms in state $s$.
    		\State Update remaining budget $B_t \leftarrow B_t - \min \{ B_t, Z_t(s)\}$.
		\EndFor \label{l:12}
        \For {state s in $C^{0}_t$}\label{l:13}
            \State Plan to pull (at least) $X_t^N(s, 1) \leftarrow \min\{B_t, Z_t(s), \lfloor N x_t(s, 1) \rfloor\}$ arms in state $s$.
            \State Store the number of undecided arms $U_t^N(s) \leftarrow Z_t^N(s) - X_t^N(s, 1)$.
            \State Update the remaining budget $B_t \leftarrow B_t - \min\{B_t, Z_t(s), \lfloor N x_t(s, 1) \rfloor \}$.
        \EndFor\label{l:16}
        \For {state s in $C^{0}_t$}\label{l:14}
    		\State Plan to pull $\min\{B_t, U_t(s)\}$ additional undecided arms in state $s$.
    		\State Update $X_t^N(s, 1) \leftarrow X_t^N(s, 1) + \min\{B_t, U_t(s)\}$.
    		\State Update $B_t \leftarrow B_t - \min \{ B_t, U_t(s)\}$.
		\EndFor \label{l:17}
        \For {state s in $C^{-}_t$}\label{l:18}
    		\State Plan to pull $X_t^N(s, 1) \leftarrow \min\{B_t, Z_t(s)\}$ arms out of $Z_t(s)$ arms in state $s$.
    		\State Update remaining budget $B_t \leftarrow B_t - \min \{ B_t, Z_t(s)\}$.
		\EndFor \label{l:21}
		\State For each $s$, pull $X_t^N(s,1)$ arms in state $s$ (as planned above)
    \EndFor
	\end{algorithmic}
\end{algorithm}
\newpage

With this definition in place, we now state the main result of this section: that fluid-priority policies are diffusion regular, implying they have an $O(\sqrt{N})$ optimality gap by Theorem \ref{diffusion_optimal}.
\begin{theorem}\label{th:3}
    Any fluid-priority policy $\pi$ is diffusion regular and its optimality gap is $O(\sqrt{N})$.
\end{theorem}

The proof of Theorem \ref{th:3} is in the Appendix.

\section{Non-degeneracy Condition: Achieving an O(1) Optimality Gap}\label{se-7}
This section presents our main contribution: that fluid-priority policies achieve an $O(1)$ optimality gap under a non-degeneracy condition. We first define and discuss this condition and then show this result.

To motivate this non-degeneracy condition, consider a fluid-priority policy and another policy motivated by the relaxed problem (\ref{KP}) in which the  
the almost-sure budget constraint ($|\textbf{a}_t| = \alpha_t N$) has been relaxed.
This so-called ``budget-relaxed'' policy first categorizes states into fluid-active, fluid-neutral, and fluid-inactive categories in the same way as its corresponding fluid-priority policy.
It pulls all fluid-active arms (even if this would exceed the budget). 
If budget remains, it then pulls fluid-neutral arms in the same way as its corresponding fluid-priority policy.
It does not pull any fluid-inactive arms, even if budget remains after fluid-active and fluid-neutral arms are pulled. 

Pulling all fluid-active arms and idling fluid-inactive arms 
is exactly the property required for any feasible policy to be optimal in the LP relaxation (\ref{KP}).
Thus, this budget-relaxed policy's reward is 
close to the relaxed problem's optimal reward (Lemma \ref{le-6}).
Moreover, it behaves identically to its corresponding fluid-priority policy (Lemma \ref{le-5}) except on a specific ``budget violation'' event: that the number of fluid-active arms exceeds the budget, or the number of fluid-active and fluid-neutral arms together fail to exceed the budget.
The probability of budget-violation allows us to bound the optimality gap for fluid-priority policies by comparing them with their budget-relaxed versions.

The non-degeneracy condition (Definition \ref{non-deg}) characterizes the probability of budget violation: when it is met, the expected number of fluid-active arms is {\it strictly} below the budget and the expected number of fluid-active and fluid-neutral arms is {\it strictly} above the budget. Thus, using concentration bounds, 
problems meeting the non-degeneracy condition are ones in which the probability of budget violation vanishes 
exponentially fast as $N$ grows (Lemma \ref{tail}). 
As a result, in such problems, the fluid-priority policy behaves the same as its budget-relaxed version with high probability for large $N$. We use this fact to show an $O(1)$ optimality gap in Theorem \ref{thm-4}.

The rest of this section is organized as follow: we first formally introduce budget-relaxed policies, then define the non-degeneracy condition, and finally prove fluid-priority policies achieve an $O(1)$ optimality gap when this condition holds. In addition, we show that no index policy can strictly outperform all fluid-priority policies. At the end of the section, we propose an algorithm to verify whether non-degeneracy condition holds, and search for a non-degenerate occupation measure when exists.

\subsection{Budget-relaxed fluid-priority policies}
Given a fluid-priority policy $\pi_F$, its budget-relaxed version $\pi_R$ is defined formally by Algorithm \ref{al:2}. Similar to $\pi_F$, $\pi_R$ first classifies states into three categories: fluid-active, fluid-neutral and fluid-inactive, using the same occupation measure as $\pi_F$. Then, $\pi_R$ sorts states in each category in order of decreasing priority-score (line \ref{l:3}), using the same priority score as $\pi_F$. Afterwards, $\pi_R$ pulls all arms in the fluid-active category (lines \ref{2l:9} - \ref{2l:12}), exceeding the budget if necessary. If there is still budget remaining, $\pi_R$ iterates over each state $s$ in the fluid-neutral category $C^{0}_t$. It pulls arms in this state until the number pulled reaches the quantity $\lfloor N x_t(s,1) \rfloor$ derived from the optimal occupation measure, no arms remain in this state, or we reach the budget. Unpulled arms in each such state are called ``undecided''. (lines \ref{2l:13} - \ref{2l:16}).
Finally, $\pi_R$ iterates over each state in the fluid-neutral category $C^{0}_t$ again and pulls undecided arms until either the budget is met or all undecided arms are pulled (line \ref{2l:14} - \ref{2l:17}). Notice $\pi_R$ idles all arms in the fluid-inactive category, even if budget remains.

\begin{algorithm}
	\caption{Budget-relaxed fluid-priority policies}\label{al:2}
	\hspace*{\algorithmicindent} \textbf{Input:} optimal occupation measure $(x_t(s, a))_{t \in [T], s \in S, a \in A}$, priority-score function $\{ \mathcal{P}_t \}_{t \in [T]}$.
	\begin{algorithmic}[1]
		\For {$t = 1, 2, ..., T$}
	    \State Classify states into fluid-active ($C^+_t$), fluid-neutral ($C^0_t$) and fluid-inactive ($C^-_t$) categories based on the occupation measure, according to equation (\ref{category}).
		\State Observe there are $Z_t(s)$ arms in state $s$, and remaining budget $B_t = \lfloor \alpha_t N \rfloor$.
		\State Each of the for loops below iterates over states in decreasing order of $\mathcal{P}_t(s)$.\label{l:3}
		
		\For {state s in $C^{+}_t$}\label{2l:9}
    		\State Plan to pull $X_t^N(s, 1) \leftarrow Z_t(s)$ arms out of $Z_t(s)$ arms in state $s$.
    		\State Update the remaining budget $B_t \leftarrow B_t - Z_t(s)$.
		\EndFor \label{2l:12}
		\If {$B_t \leq 0$}
		    \State continue
		\EndIf
        \For {state s in $C^{0}_t$}\label{2l:13}
            \State Plan to pull (at least) $X_t^N(s, 1) \leftarrow \min\{B_t, Z_t(s), \lfloor N x_t(s, 1) \rfloor \}$ in state $s$.
            \State Store the number of undecided arms $U_t^N(s) \leftarrow Z_t^N(s) - X_t^N(s, 1)$.
            \State Update the remaining budget $B_t \leftarrow B_t - \min\{B_t, Z^N_t(s), N x_t(s, 1)\}$.
        \EndFor\label{2l:16}
        \For {state s in $C^{0}_t$}\label{2l:14}
    		\State Plan to pull $\min\{B_t, U_t(s)\}$ additional undecided arms in state $s$.
    		\State Update $X_t^N(s, 1) \leftarrow X_t^N(s, 1) + \min\{B_t, U_t(s)\}$, 
    		\State Update $B_t \leftarrow B_t - \min \{ B_t, U_t(s)\}$.
		\EndFor \label{2l:17}
	\State For each $s \in C_t^+ \cup C_t^0$, pull $X_t^N(s, 1)$ arms in state $s$ (as planned above). Idle all arms in $C_t^-$.
    \EndFor
	\end{algorithmic}
\end{algorithm}
\newpage

Policy $\pi_R$ behaves the same as the corresponding fluid-priority policy $\pi_F$ with the same occupation measure and priority-score function, except $\pi_R$ pulls all arms in the fluid-active category and idles all arms in the fluid-inactive category regardless of the budget constraint. We state this formally in the following lemma, whose proof is in the Appendix.
\begin{lemma}\label{le-5}
    Define event $\Delta_t$ as
    \begin{align*}
        \Delta_t := \Big\{ \sum_{s \in C^{+}_t} Z^N_t(s) \leq \alpha_t N \leq \sum_{s \in C^{+}_t} Z^N_t(s) + \sum_{s \in C^{0}_t} Z^N_t(s) \Big\}.
    \end{align*}
    Then $\pi_R(t, Z_t^N) = \pi_F(t, Z_t^N)$ on the event $\Delta_t$.
\end{lemma}
We write the complement of $\Delta_t$ as $\Delta_t^c$
and refer to this as a ``budget violation'' event. 


\subsection{Non-degeneracy} The non-degeneracy condition states that the fluid-neutral category is not empty, which is sufficient to prove the budget-violation events, $\Delta_t^c$, are probabilistically negligible. 
Roughly speaking, non-emptiness of $C^{0}_t$ guarantees that the occupation measure satisfies,
\begin{align*}
    \sum_{s \in C^{+}_t}{z_t(s)} < \alpha_t < \sum_{s \in C^{+}_t}{z_t(s)} + \sum_{s \in C^{0}_t}{z_t(s)}.
\end{align*}
Along with both the budget-relaxed fluid-priority policy and the fluid-priority policy being fluid consistent, the number of arms in state $s$, $Z_t^N(s)$, is roughly proportional to $z_t(s)$, with excursions described by a central limit theorem. Thus the probability of event $\Delta_t$ approaches $1$ exponentially fast with $N$ by concentration inequalities. 

We formally state this in the following definition, first defining non-degeneracy, and then stating in Lemma~\ref{tail} that non-degeneracy implies that budget violations are probabilistically negligible for large $N$. 
The proof of Lemma \ref{tail} is in the Appendix.

\begin{definition}\label{non-deg}
We say an optimal occupation measure $(x_t(s, a))_{t \in [T], s \in S, a \in A}$ is {\it non-degenerate} if 
\begin{align*}
    |C^{0}_t| \geq 1, \forall t \in [T].
\end{align*}
Otherwise, we call it {\it degenerate}. We also call a fluid-priority policy {\it non-degenerate} ({\it degenerate}) when its associated occupation measure is non-degenerate (degenerate).
\end{definition}

\begin{lemma}\label{tail}
    If an optimal occupation measure $(x_t(s, a))_{t \in [T], s \in S, a \in A}$ is non-degenerate, then for any priority-score functions $\{ \mathcal{P}_t \}_{t \in [T]}$ and the corresponding fluid-priority policy $\pi_F$ and budget-relaxed policy $\pi_R$,
    there exists a constant $\delta > 0$ and a constant $L$ such that 
    \begin{align*}
        \max\{\mathbb{P}_{\pi_R}(\Delta_t^c),
        \mathbb{P}_{\pi_F}(\Delta_t^c)\}
        \leq L \exp(-\delta N),
    \end{align*}
   for all $t \in [T]$ and all $N$. 
\end{lemma}

Empirically, one can check the non-degeneracy condition for a given optimal occupation measure $x^*$ returned by solving the linear programming relaxation \eqref{LP} with a commercial LP solver. Recalling from \eqref{category} that states $s$ in the fluid-neutral category are those with both $x_t(s,1)>0$ and $x_t(s,0)>0$, we check whether $x^*$ is non-degenerate by assessing whether there is at least one such state for each $t$.



\subsection{Main result} 
We now state and prove this section's main result: 
a fluid-priority policy achieves an $O(1)$ optimality gap when it is non-degenerate.  Before that, we need one last building block: the budget-relaxed policies' reward deviates from the relaxed problem's optimal reward by at most $O(1)$ under non-degeneracy.

We show this in the following lemma, whose proof is in the appendix. There are two main ideas in the proof.
First, recall that the budget-relaxed fluid-priority policy $\pi_R$ pulls all arms in $C^{+}_t$ and no arms in $C^{-}_t$. Thus, its decisions are optimal under a Lagrangian relaxation of Problem~\ref{KP} in which the budget constraint on the expected number of arms pulled is replaced by a well-chosen linear penalty (in this Lagrangian relaxation, fluid-neutral arms does not affect optimality as the incremental reward is offset by the linear penalty).
Second, the fact that $\pi_R$ pulls a number of arms equal to the almost-sure budget constraint with high probability ensures that it nearly satisfies the constraint in Problem~\ref{KP} on the expected budget. This fact causes the linear penalty to be nearly 0. Finally, the fact that the Lagrangian relaxation is the sum of the unpenalized reward, which we call $V_N(\pi_R)$, and this penalty imply that the $V_N(\pi_R)$ is within a constant of the value of Problem~\ref{KP}, $\hat{V}_N^*$.

\begin{lemma}\label{le-6}
Let $V_N(\pi_R) =\max_{\pi} \mathbb{E}_{\pi} \sum_{t=1}^T R_t(\mathbf{s}_t, \mathbf{a}_t)$ for a budget-relaxed fluid priority policy $\pi_R$.
If an optimal occupation measure $(x_t(s, a))_{t \in [T], s \in S, a \in A}$ is non-degenerate, then for any priority-score functions $\{ \mathcal{P}_t \}_{t \in [T]}$, the corresponding budget-relaxed fluid-priority policy $\pi_R$ satisfies
\begin{align*}
    |\hat{V}_N^* - V_N(\pi_R)| \leq m,
\end{align*}
where $m$ is a constant not depending on $N$.
\end{lemma}

Now we are ready to state and prove our main result: that a fluid-priority policy achieves an $O(1)$ optimality gap when it is non-degenerate.
A fluid-priority policy $\pi_F$'s optimality gap can be bounded by first comparing the reward $V_N(\pi_F)$ with the reward of its corresponding budget-relaxed policy $\pi_R$. Combining the fact that $\pi_F$ deviates from $\pi_R$ with negligible probability (Lemma \ref{tail}) and that $\pi_R$'s reward deviates $O(1)$ from $\hat{V}_N^*$ (Lemma \ref{le-6}), $V_N(\pi_F)$ is at most $O(1)$ away from $\hat{V}_N^*$. 

\begin{theorem}\label{thm-4}
    If an optimal occupation measure $(x_t(s, a))_{t \in [T], s \in S, a \in A}$ is non-degenerate, then for any priority-score functions $\{ \mathcal{P}_t \}_{t \in [T]}$, the corresponding fluid-priority policy $\pi_F$ satisfies
    \begin{align*}
        \hat{V}_N^* - V_N(\pi_F) \leq m,
    \end{align*}
    where $m$ is a constant not depending on $N$.
\end{theorem}

\proof{Proof of Theorem \ref{thm-4}}


Under $\pi_F$, the reward is
\begin{align*}
    V_N(\pi_F) = \mathbb{E}_{\pi_F}\sum_{t = 1}^T \sum_{s \in S, a \in A} r_t(s, a) X_t^N(s, a).
\end{align*}

Under $\pi_R$, the reward is
\begin{align*}
    V_N(\pi_R) &= \mathbb{E}_{\pi_R}\sum_{t = 1}^T \sum_{s \in S, a \in A} r_t(s, a) X_t^N(s, a).
\end{align*}

Denote $\Omega_T := \Delta_1 \cap \Delta_2 \cap ... \cap \Delta_T$.  On this event, $\pi_R$ and $\pi_F$ produce identical decisions by Lemma~\ref{le-5}. Using this in the second line below, we have:
\begin{align*}
    V_N(\pi_R) - V_N(\pi_F)
    &= \mathbb{E}_{\pi_R} \Big[ 1_{\Omega_T} \sum_{t = 1}^T \sum_{s \in S, a \in A} r_t(s, a) X_t^N(s, a) \Big] + \mathbb{E}_{\pi_R} \Big[ 1_{\Omega^c_T}\sum_{t = 1}^T \sum_{s \in S, a \in A} r_t(s, a) X_t^N(s, a) \Big] \\
    & - \mathbb{E}_{\pi_F} \Big[ 1_{\Omega_T} \sum_{t = 1}^T \sum_{s \in S, a \in A} r_t(s, a) X_t^N(s, a) \Big]  - \mathbb{E}_{\pi_F} \Big[ 1_{\Omega^c_T} \sum_{t = 1}^T \sum_{s \in S, a \in A} r_t(s, a) X_t^N(s, a) \Big] \\
    & = \mathbb{E}_{\pi_R} \Big[ 1_{\Omega^c_T} \sum_{t = 1}^T \sum_{s \in S, a \in A} r_t(s, a) X_t^N(s, a) \Big] - \mathbb{E}_{\pi_F} \Big[ 1_{\Omega^c_T} \sum_{t = 1}^T \sum_{s \in S, a \in A} r_t(s, a) X_t^N(s, a) \Big] \\
    & \leq \mathbb{E}_{\pi_R} \Big[ 1_{\Omega^c_T} \sum_{t = 1}^T \sum_{s \in S, a \in A} | r_t(s, a)| |X_t^N(s, a) | \Big]  + \mathbb{E}_{\pi_F}  \Big[ 1_{\Omega^c_T} \sum_{t = 1}^T \sum_{s \in S, a \in A} | r_t(s, a)| |X_t^N(s, a) | \Big].
\end{align*}

Inequalities $0 \leq X^N_t(s) \leq N$ and $0 \leq Y^N_t(s) \leq N$ then imply
\begin{align*}
    V_N(\pi_R) - V_N(\pi_F)
    & \leq \mathbb{E}_{\pi_R} \Big[ 1_{\Omega^c_T} \sum_{t = 1}^T \sum_{s \in S, a \in A} | r_t(s, a)| N \Big] + \mathbb{E}_{\pi_F} \Big[ 1_{\Omega^c_T} \sum_{t = 1}^T \sum_{s \in S, a \in A} | r_t(s, a)| N \Big] \\
    & \leq 2\left(\mathbb{E}_{\pi_R} [1_{\Omega^c_T}]  +  \mathbb{E}_{\pi_F} [1_{\Omega^c_T}]\right) T |S| r^* N,
\end{align*}
where $r^* := \max_{t, s, a} |r_t(s, a)|$.

Then, applying Lemma \ref{tail} and $\mathbb{P}_{\pi}[\Omega_T^c] \le \sum_{t = 1}^T \mathbb{P}_{\pi}[\Delta_t^c]$, we have: 
\begin{align*}
    V_N(\pi_R) - V_N(\pi_F)
    & \leq 2T |S| r^* N \sum_{t = 1}^T \mathbb{P}_{\pi_R} [\Delta^c_t] + \mathbb{P}_{\pi_F}[\Delta^c_t] \\
    & \leq 4 T^2 |S| r^* N L \exp(-\delta N).
\end{align*}
Finally, applying Lemma \ref{le-6} concludes the proof.
\hfill$\Box$


\subsection{The best fluid-priority policy is at least as good as the best index policy} 
Here we compare fluid-priority policies against index policies \citep{whittle1980multi, gittins2011multi}.
An index policy assigns each state an ``index'' and prioritizes each arm based on the index of its current state from high to low, pulling arms until we exhaust the current period's budget.  

A policy can be both a fluid-priority policy and index policy. This occurs if there is at most one fluid-neutral state in any period and indices of all fluid-active states are higher than those of all fluid-neutral states, which are higher in turn than the indices of all fluid-inactive states. There are, however, index policies that are not fluid priority policies, and vice versa.  If the indices do not respect the ordering implied by the fluid-active, fluid-neutral, and fluid-inactive categories then this index policy is not a fluid priority policy. Also, if multiple fluid-neutral states can be occupied in one period, a fluid-priority policy will allocate pulls across these arms in accordance with an occupation measure and in a way that is different from the strict prioritization used by an index policy.

Since index policies are widely known and used, it is instructive to compare them with fluid-priority policies.
The discussion above shows that the difference between $\hat{V}_N^*$ (the optimal objective of the relaxation) and the value of a fluid-priority policy $V_N(\pi_{F})$ is bounded above by a constant when $\pi_F$ is non-degenerate, i.e., that 
$\sup_N \hat{V}_N^* - V_N(\pi_F)$ is finite.
The following proposition shows that the best fluid priority is at least as good as the best index policy, when measured by
$\sup_N \hat{V}_N^* - V_N(\pi) \in \mathbb{R} \cup \{\infty\}$, regardless of whether non-degeneracy holds.

\begin{proposition}
\label{prop-1}
Consider an index policy $\pi_I$ such that 
$\sup_N \hat{V}_N^* - V_N(\pi_I) < \infty$.
Then, there exists a fluid priority policy $\pi_F$ such that 
$$\sup_N \hat{V}_N^* - V_N(\pi_F) 
\le \sup_N \hat{V}_N^* - V_N(\pi_I).$$
\end{proposition}
The proof of Proposition \ref{prop-1} is in the Appendix.

\subsection{Choice of Occupation Measure}\label{select-occupation-measure}

The above discussion of fluid-priority policies and degeneracy applies to any optimal occupation measure. Multiple optimal occupation measures may exist, some degenerate and others not. In this situation, a fluid-priority policy constructed using a non-degenerate optimal occupation measure is guaranteed to have an $O(1)$ optimality gap while another constructed using a degenerate one is not. A natural question then arises: how can we determine whether a non-degenerate optimal occupation measure exists and how can we select one if it does? Here we describe a computational procedure for answering this question.

First, observe from \eqref{LP} that any convex combination of optimal occupation measures is also optimal.
Thus, suppose we can find a collection of optimal occupation measures, $x^{*,k}$, $k \in [K]$, such that, for each $t$, there is either (1) a state $s$ that is fluid-neutral under some $k$, or (2) there is a state $s$ that is fluid-active under some $k$ and fluid-inactive under another $k$.
Then any convex combination with strictly positive weight on each $k$ is non-degenerate. 
We describe an algorithm for finding such a collection, if it exists, or establishing that it does not. 

To accomplish this, first solve the LP \eqref{LP}, call the solution $x^{*,1}$, and record its optimal value for later use. Assess for each $t$ whether there is a state $s$ satisfying $x^{*}_t(s,0) >0$ and $x^{*}_t(s,1) > 0$.
If all $t$ satisfy this condition, then we have found a non-degenerate optimal occupation measure.

Otherwise, we will continue iteratively in our search.
In each stage $k$, 
we will maintain a collection of solutions $\{x^{*,k'} : k' = 1,\ldots, k\}$ and a set of times $A_k \subseteq [T]$.
$A_k$ contains those times for which we have not yet been able to construct a fluid-neutral state. 
Formally, a time $t$ is in $A_k$ if and only one of the following holds: (1) all states are fluid-active at $t$ in all $x^{*,k'}$, $k'\le k$; or (2) all states are fluid-inactive at $t$ in all $x^{*,k'}$, $k'\le k$. 
If $A_k$ is empty, then a non-degenerate optimal occupation measure can be constructed as a convex combination of $\{x^{*,k'} : k' = 1,\ldots, k\}$ using 
strictly positive weights on every solution in this collection.
If $A_k$ is not empty, we will then attempt to construct an optimal occupation measure that, when added to our collection of solutions, causes $A_{k+1}$ to be a strict subset of $A_k$.

Toward this goal, in stage $k$, choose $t_k \in A_k$. This will be the time that we seek to remove from $A_k$ in constructing $A_{k+1}$.
Let $C^{+,k}$ contain all of the states for which $x^{*,k'}_{t_k}(s,1) > 0$ and $x^{*,k'}_{t_k}(s,0) = 0$    for all $k' \le k$. These are the states that are fluid-active at time $t_k$ for all previously computed optimal occupation measures.
Then solve a linear program minimizing  
$\sum_{s \in C^{+,k}} x_{t_k}(s,1)$
subject to 
all of the constraints in \eqref{LP}
and the linear constraint that the objective in \eqref{LP} is equal to its optimal value recorded above.
Call the solution $x^{*,k+1}$.

This linear program assesses whether there is an optimal occupation measure $x^{*,k+1}$ satisfying 
$\sum_{s \in C^{+,k}} x_{t_k}(s,1) < \alpha_k$.
If no such $x^{*,k+1}$ exists, then this establishes that all optimal occupation measures are degenerate.
Otherwise, if we find such a $x^{*,k+1}$, then we add it to our collection of solutions. We also construct $A_{k+1}$ by removing the time $t_k$ from $A_k$. We additionally remove any other times $t$ for which the new solution $x^{*,k+1}$ provides a state whose category at that time $t$ is different from those in the previous solutions $x^{*,k'}$, $k'\le k$.

If $A_{k+1}$ is the empty set, then this implies that there is a non-degenerate optimal occupation measure. We set $K=k+1$ and construct it as described above from the collection $\{x^{*,k'} : k \le K\}$.

\section{Numerical Experiments}\label{numeric}

This section evaluates the performance of fluid-priority and other policies on three problems, leveraging both simulation experiments, computational investigations and our earlier theoretical results.

\S\ref{numeric-1} studies a classical problem: the finite-horizon Bayesian bandit with Bernoulli rewards. It compares fluid-priority policies against the widely used UCB and Thompson Sampling policies. Fluid-priority policies are first shown to substantially outperform both methods in numerical experiments. We then show that the optimality gap for UCB and Thompson sampling is $\Omega(N)$, while it is $O(1)$ for the fluid-priority policy evaluated as it is non-degenerate in this problem.

\S\ref{numeric-2} considers an active learning problem based on \cite{chen2013optimistic} in which an algorithm allocates crowd workers (e.g., participants on Amazon's Mechanical Turk) to image labeling tasks to support learning an accurate classifier. We show via numerical experiments that fluid-priority policies outperform a previously proposed state-of-the-art policy (Optimistic Knowledge Gradient, \citealt{chen2013optimistic}) specifically designed for this problem.

\S\ref{numeric-3} shows via direct computation that the dynamic assortment problem previously studied in \cite{brown2020index} satisfies the non-degeneracy condition. This and our main theoretical result shows that fluid-priority policies have an $O(1)$ optimality gap, explaining the hitherto poorly understood performance of Lagrangian index policies first noticed in numerical experiments in \cite{brown2020index}.

\subsection{Bayesian bandit with Bernoulli rewards}\label{numeric-1}
This section evaluates fluid-priority policies performance on the Bayesian bandit problems, which is a standard benchmark in the bandit literature. 
While the problem is not restless, and so does not demand the full capabilities of our proposed fluid-priority policies, it allows us to study benchmarks designed for non-restless settings. Problems with restless arms are studied later.

We compare the performance of fluid-priorities against Upper Confidence Bound (UCB, \citealt{agrawal1995sample}) and Thompson Sampling (TS, \citealt{agrawal2012analysis}) policies and show that the fluid-priority policy achieves an $O(1)$ optimality gap while the optimality gaps of both UCB and TS grow linearly with the number of arms.
While UCB and TS are well-known for having a logarithmic asymptotic performance guarantee of $O(N\log(T))$, this is linear in $N$. (It also applies to a slightly different problem setting than the one we study here:
a stochastic frequentist setting with one pull period and where regret is measured with respect to the policy that pulls the best arm.)
Thus, the classical regret guarantee for these policies is not inconsistent with our finding that these policies have a $\Omega(N)$ optimality gap in a Bayesian analysis with multiple pulls per period. 

This suggests that when $T$ is small and $N$ is large, and where prior information supports the use of a Bayesian analysis, there is significant value in using fluid priority policies over UCB or TS, and in using a Bayesian finite-horizon analysis rather than a stochastic frequentist analysis.

\textbf{Problem Setup:}
There are $N$ arms in total, of which we may pull at most $\lfloor N/3 \rfloor$ in each of $T$ periods. 
Before any arms are pulled, each arm $i$ has a parameter 
$\theta_i$ sampled independently from the Bayesian prior distribution on the arm's reward. This prior distribution is uniform with support $[0, 1]$. Then, conditioning on $\theta_i$, each arm $i$'s rewards are generated when pulled as conditionally independent Bernoulli random variables with a common parameter $\theta_i$.
 Our objective is to maximize the expected total reward collected across all periods.
 
 This problem is similar to the more widely-studied stochastic bandit, except that the arm's reward is drawn at random from the prior. The expected reward calculated can be understood as the average-case reward over stochastic bandit problem instances, i.e., over ($\theta_i : i)$, where the weight on a particular instance 
($\theta_i : i)$ is proportional to its density under the prior.

The non-degeneracy condition holds in this problem for both horizons $T=15$ and $T=20$. We verified this numerically by solving the linear program \eqref{LP} and confirming that there is at least one state with a strictly positive occupation measure  in the fluid-neutral category in each period.

\textbf{Policy Implementation:}
We briefly discuss how we implement UCB, TS and fluid-priority policies in these experiments.

UCB tracks the posterior belief on $\theta_i$ for each arm $i$ based on the arm's past reward realization, and calculates an upper confidence bound for $\theta_i$ as $\mu_i + \delta \sigma_i$, where $\delta$ is a fixed parameter, $\mu_i$ is the mean of the posterior belief and $\sigma_i$ is the standard deviation. The top $\lfloor N / 3 \rfloor$ arms ranked by their upper confidence bound are selected to be pulled. We run UCB with $\delta$ varying from $0.1$ to $1$ and report results for the one with the best expected reward ($\delta^* = 0.5$) in both experiments. 

TS also tracks the posterior belief on $\theta_i$ for each arm $i$. 
At each period, TS samples a value from each arm's posterior belief on $\theta_i$, then pulls the $\lfloor N / 3 \rfloor$ arms with the highest sampled values. 

The fluid-priority policy is constructed as follows. 
First, to construct the optimal occupation measure, we solve the relaxed problem (\ref{LP}) and fetch its solution. 
Second, to construct the priority-score function,
we use a Lagrangian-relaxation approach similar to \cite{hu2017asymptotically} and \cite{brown2020index}: we solve the min-max problem
\begin{equation}\label{LR-P}
\begin{split}
    (\lambda^*_1, ..., \lambda^*_T)
    \leftarrow \min_{\lambda_1, ..., \lambda_T} \max_{\pi} \ \mathbb{E}_{\pi} \sum_{t=1}^T r_t(s_t, a_t) + \lambda_t (\alpha_t - a_t),
\end{split}
\end{equation}
where the inner $\max$ can be solved via dynamic programming and the outer $\min$ can be solved via the subgradient method.
Then we compute the $Q-$function based on the optimal Lagrangian multiplier $(\lambda^*_1, \lambda^*_2, ..., \lambda^*_T)$ iteratively:
\begin{align*}
    Q_t(s, a) = r_t(s, a) - \lambda_t a + \sum_{s'} p_{t}(s, a, s') \max_{a'} Q_{t + 1}(s', a') \text{ for $0 \le t \le T - 1$,}
\end{align*}
with $Q_T(s, a) = r_T(s, a) - \lambda_T a$, and construct the priority-score function as $\mathcal{P}_t(s) = Q_t(s, 1) - Q_t(s, 0)$.
Finally, we plug the optimal occupation measure and the score-function into Algorithm \ref{al:1} to construct the fluid-priority policy.

\textbf{Numerical Experiments:}
We compare the just-described fluid-priority policy against UCB and TS using two different time horizons $T$ of 15 and 20.
Figure 1a displays results for $T = 15$ while Figure 1b shows results for $T = 20$.
In both experiments, we iteratively double the number of arms (from $N = 300$ to $38400$) and plot an upper bound on the optimality gap. 
This bound on the optimality gap is computed by first computing the value of the relaxed problem $\hat{V}_N^*$ (which is an upper bound on the value of the optimal policy) and then subtracting the value of the UCB, TS or fluid-priority policy in question estimated via simulation.
We compare this upper bound across policies instead of the exact optimality gap because computing the exact optimality gap would require knowing the value of the optimal policy, which would take time exponential in $N$, as discussed in \S\ref{approx}. We use $50N$ replications to estimate a policy's value when there are $N$ arms. We use more samples when there are  more arms because having more arms increases the variance of a policy's reward.
We also compute a confidence interval on this upper bound, computed as the difference between $\hat{V}^*_N$ 
and the upper and lower limits of a confidence interval on the policy's expected reward.

\begin{figure}[htb]%
    \centering
    \subfloat[$T = 15$]{{\includegraphics[width=8cm, height = 6 cm]{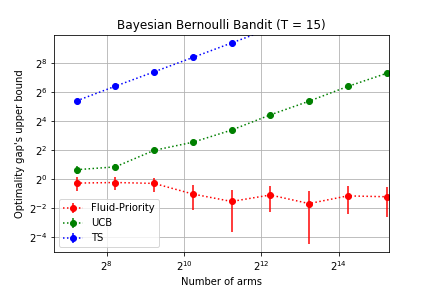}}}%
    \subfloat[$T = 20$]{{\includegraphics[width=8cm, height = 6 cm]{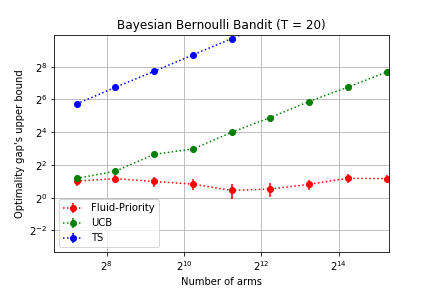}}}%
    \caption{Bayesian bandit with Bernoulli rewards. An upper bound on the optimality gap (relaxed problem's expected total reward minus a simulation-based estimate of reward) vs number of arms $N$, for the finite-horizon Bayesian multi-armed bandit with horizons $T=15$ (left) and $T=20$ (right). The fluid-priority policy has its optimality gap bounded above by a constant while UCB and Thompson sampling have optimality gaps that grow linearly with the number of arms.}
    \label{fig:vannilla_bandit}%
\end{figure}


Figure \ref{fig:vannilla_bandit} compares the performance of fluid-priority, UCB and TS policies. For both time horizons $T$ of 15 and 20, the fluid-priority policy performs significantly better than UCB and TS, especially for large $N$.
The fluid-priority policy's reward differs from the optimal policy's reward by at most $1$ for $T=15$ and at most $2$ for $T=20$ even when there are 38400 arms available. 
UCB outperforms TS in both time horizons, perhaps due to the tuning of UCB's hyperparameter.

These results are consistent with Theorem \ref{thm-4} and our numerical validation that the non-degeneracy condition is satisfied, which implies 
that the fluid-priority policy's optimality gap is bounded across all values of $N$. 
In contrast, the optimality gap for both UCB and Thompson sampling grows with $N$ as shown in Figure \ref{fig:vannilla_bandit}.

Proposition~\ref{ucb_ts_regret} provides an additional
analysis to confirm theoretically that UCB and Thompson sampling have optimality gaps that grow linearly in $N$.
The proof of this proposition, which may be found in the appendix, defines an iterative algorithm over $t$ to calculate the occupation measure for UCB and TS in the large $N$ limit. We then use this algorithm to compute this occupation measure for specific values for $T$ and compare it to the optimal occupation measure. We find that the occupation measures are suboptimal for the values of $T$ used in these experiments, implying that UCB and TS are not fluid-consistent and their optimality gaps are $\Omega(N)$. These values of $T$ are representative, and UCB and TS have $\Omega(N)$ optimality gaps for other $T$ as well.
\begin{proposition}\label{ucb_ts_regret}
The optimality gap for both UCB and TS is $\Omega(N)$ for $T = 15$ and $T = 20$.
\end{proposition}

\subsection{Crowdsourced Labeling}\label{numeric-2}
This section evaluates a fluid-priority policy's performance on an active learning problem introduced by \cite{chen2013optimistic} focused on the allocation of crowd workers for accurate image classification. We compare its performance against the Optimistic Knowledge-Gradient \citep{chen2013optimistic}, a method specifically designed for this problem, and the Online Knowledge-Gradient \citep{ryzhov2012knowledge}. The fluid-priority policy outperforms both methods significantly.

We formulate the crowdsourced labeling problem as follows.
Suppose there are $N$ images needing binary labels (e.g., whether this is a picture of a pedestrain or not) to support training of an automatic image classifier that will be built later. We ask crowd workers to label these images.
This approach to ``jump starting'' machine learning classifiers with labels from crowd workers is common in indutry \citep{chen2013optimistic}.
Each image $i$ has a true underlying binary class, along with an associated probability $p_i$ that a crowd worker will label the image with the correct class.
A crowd worker may provide an incorrect label because, e.g., the image is blurry or the worker is distracted.
We assume $p_i > 1/2$ (following \citealt{chen2013optimistic}), i.e. the majority of crowd workers give the correct label.
We use an independent prior belief for each image's $p_i \sim U[1/2, 1]$.
We are allowed to request $T=7$ batches of labels from crowd workers, with up to $\lfloor N / 4 \rfloor$ images per batch. After the last batch, we estimate each image's class via majority vote, which is also the class with maximum probability under posterior.

\begin{figure}[H]
	\centering
	\subfloat[]{\includegraphics[width=8cm, height = 6 cm]{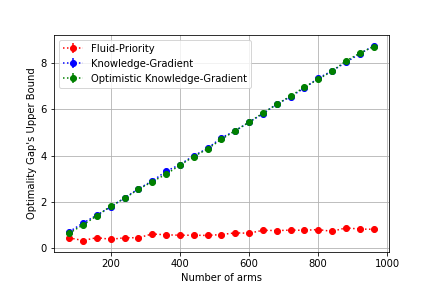}\label{fig:f1}}
	\hfill
	\subfloat[]{\includegraphics[width=8cm, height = 6 cm]{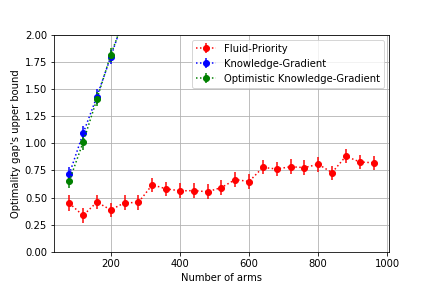}\label{fig:f2}}
	\caption{Crowdsourced labeling. 
An upper bound on the optimality gap (relaxed problem's expected total reward minus a simulation-based estimate of reward) vs. number of arms $N$. The left and right panel show the same data but using different scales for the y-axis.
Both Knowledge Gradient and Optimistic Knowledge Gradient have optimality gaps that seem to grow linearly. The fluid-priority policy has an optimality gap that is $O(\sqrt{N})$ because the non-degeneracy condition does not hold in this problem.}

	\label{crowdsource}
\end{figure}

Figure \ref{crowdsource}
compares the fluid-priority policy against the Online Knowledge-Gradient and Optimistic Knowledge-Gradient methods as we vary the number of arms $N$, reporting an upper bound on the optimality gap for each policy computed in the same way as Section~\ref{numeric-1}. 
The fluid-priority policy seems to perform extremely well, and incorrectly classifies at most 1 more image on average than the optimal policy even when there are 1000 images' labels to be learned. Online Knowledge-Gradient and Optimistic Knowledge-Gradient perform similarly in our experiment, and they both underperform the fluid-priority policy by wrongly classifying at least 8 more images on average with 1000 images' labels to be learned.
Even though the optimistic knowledge-gradient was designed specifically for this problem, the fluid-priority policy nevertheless has a significantly smaller optimality gap.

Figure \ref{crowdsource} is consistent with our theoretical results.
We can verify the non-degeneracy condition does not hold in this example by implementing the algorithm in \S\ref{select-occupation-measure}.
The lack of non-degeneracy implies that the optimality gap of the fluid-priority policy is $O(\sqrt{N})$. Its performance in the plot is consistent with this scaling. 
The Online Knowledge Gradient and the Optimistic Knowledge Gradient, however, seem to have suboptimality that scales linearly with $N$.

\subsection{Dynamic Assortment Optimization}\label{numeric-3}
This section discusses a dynamic assortment optimization problem proposed in Section 6.2.1 in \cite{brown2020index}. 
In this problem, \cite{brown2020index} observes empirically that the optimality gap of the policy proposed in that paper (a so-called Lagrangian policy), shown there to be $O(\sqrt{N})$, seems to stay constant with $N$, suggesting that the $O(\sqrt{N})$ bound is loose.
We first describe the problem setting and then confirm that our theoretical results provide the tighter bound suggested by these empirical results.

A retailer repeatedly chooses products to display in a selling season. The retailer has $N$ products 
but a shelf-space constraint allows
only showing $\lfloor N / 4 \rfloor$ of them in each time period. Each product, if sold, generates profit of \$1. 
The demand rate for each product $i$ is unknown to the retailer but follows a Poisson process with intensity $\gamma_i$.
The retailer holds a Bayesian prior belief on $\gamma_i$, which is Gamma-distributed with shape parameter $m_i$ and inverse scale parameter $a_i$, $\gamma_i \sim \mathrm{Gamma}(m_i,a_i)$. All products share the same prior belief $(m_i, a_i) = (1, 0.1)$.
The retailer updates these prior beliefs after observing demand realizations for displayed products using Bayes rule. 

The Gamma distribution is a conjugate prior distribution when we have Poisson-distributed observations, which causes the posterior to remain Gamma-distributed. 
More specifically, the posterior on $\gamma_i$ in time period $t$ is $\mathrm{Gamma}(m_{t,i}, a_{t,i})$ where $m_{t,i}$ and $a_{t,i}$ can be computed recursively. For a product $i$ that was displayed in time period $t$, letting $x_{t,i}$ be the demand for the product in the period,
$m_{t+1,i} = m_{t,i} + x_{t,i}$ and $a_{t+1,i}=a_{t,i}+1$.
For a product $i$ that was not displayed in $t$, 
$m_{t+1,i} = m_{t,i}$ and $a_{t+1,i}=a_{t,i}$.
At $t=0$,
$m_{0,i} = 1$ and $a_{0,i} = 0.1$.

The retailer's objective is to adaptively choose which products to display in each period subject to the shelf-space constraint to maximize the expected total profit over a selling season lasting $T$ periods. 
This is formulated as a restless bandit with time horizon $T$ where each product $i$ is an arm whose state at time $t$ is $(m_{t,i}, a_{t,i})$. 
The optimal policy and good approximate policies must balance exploration and exploitation by showing products that observed sales and the prior suggest have large $\gamma_i$ (exploitation) and also showing those for which we have little observed sales data to support learning $\gamma_i$  (exploration). This must be done cognizant of the time horizon $T$: as the remaining time in the selling season shrinks, exploration becomes unimportant. This must also be done while leveraging the prior, especially when $T$ is small and observed sales alone leave substantial uncertainty.

\cite{brown2020index}  study performance of their proposed Lagrangian policy when $T = 8$. They find their policy ``perform(s) very well for large $N$'', and produces profit ``within \$6 of the optimal value!'' when $N=16,384$.
They do not, however, offer an explanation for why the performance would be so good for a policy with an $O(\sqrt{N})$ optimality gap, the tightest bound known at the time.

Our results explain this phenomenon.
First, by solving the linear programming relaxation \eqref{LP} for this problem, we confirm that the set of fluid-neutral states is non-empty in each time period, thus confirming that the problem is non-degenerate.
Moreover, there is exactly one state in each period's fluid-priority category.
This is also observed by \cite{brown2020index}, as they mention that ``there are no scenarios where products in different states have the same priority indices''. Thus, for this optimal occupation measure, fluid-priority policies are also index policies and the Lagrangian policy is one specific example. This explains why the Lagrangian policy achieves an $O(1)$ optimality gap.

\section{Conclusion}\label{conclusion}
In this work,
we have developed fluid-priority policies, 
a class of new policies with strong theoretical guarantees and numerical performance for Bayesian restless bandits with many arms. In the limit as the number of arms $N$ grows large, fluid-priority policies always achieve an $O(\sqrt{N})$ optimality gap, matching the best existing scaling in the past literature. When the non-degeneracy condition holds, surprisingly fluid-priority policies achieve an $O(1)$ optimality gap.

Although our analysis is specific to Bayesian restless bandits , our general approach may support analysis of policies based on fluid approximations in other areas.
Fluid approximations and policies based on them are a widely used tool 
in revenue management \citep{dai2019network}, inventory control \citep{kunnumkal2011linear} and other areas. 
They are used, in particular, in Weakly Coupled Markov Decision Processes, which generalize restless bandits by allowing multiple actions for each arm and multiple resource constraints. As our analysis of Bayesian restless bandits only leverages a fluid-approximation method and a concentration property in a many-arm regime, our results are likely generalizable to problem formulations from these other areas under a similar asymptotic regime.
We also believe it is possible to extend our work to infinite-horizon problems, and to   incorporate correlated randomness across resource constraints, e.g., through resource constraint bounds given not by deterministic values but by a Markov process that is common across arms. 
In summary, we feel that the technique demonstrated here of concentration-inequality analysis of deviations from fluid policies is a powerful technique that can be broadly applied.

\newpage
\input{appendix.tex}

\bibliographystyle{informs2014}
\bibliography{References}
\end{document}

%% file: Appendix.tex
\section{Appendix}

This section provides all technical proofs not included in the main paper.
\subsection{Proof for Lemma \ref{fenchel}}
In the original formulation of the restless bandit, problem \eqref{OP}, the budget constraint $|\mathbf{a}_t| = \lfloor \alpha_t N \rfloor$ applies on each sample path. 
The relaxed problem \eqref{KP} is identical except that this constraint is replaced by the weaker one, $\mathbb{E}|\mathbf{a}_t| = \alpha_t N$.
Recalling our assumption here that $\alpha_t N$ is an integer, the right-hand sides of these two constraints are the same. (Generalizations to non-integer $\alpha_t N$ are discussed in Appendix~\ref{al:2-1}).
Thus, the set of feasible policies in \eqref{OP} is a subset of those in \eqref{KP}, implying that the value of \eqref{OP} is bounded above by that of \eqref{KP}, i.e.,
\begin{align}
    V_N^* \leq \hat{V}_N^{*}.
\end{align}
To prove $\hat{V}^*(N) = N \hat{V}^*(1)$, we use a Lagrangian Relaxation similar to \cite{farias2011irrevocable, guha2008sequential} as the key idea in the following argument.

Through straightforward imitation of the proof of the Fenchel Duality Theorem  \citep{rockafellar1970convex},
\begin{align}\label{eq-8}
    \max_{\pi} \min_{\lambda_{1: T}} \mathbb{E}_{\pi} \sum_{t = 1}^T R_t(\mathbf{s}_t,\mathbf{a}_t) + \lambda_t (\alpha_t N - |\mathbf{a}_t|) = \min_{\lambda_{1: T}} \max_{\pi} \mathbb{E}_{\pi} \sum_{t = 1}^T R_t(\mathbf{s}_t,\mathbf{a}_t) + \lambda_t (\alpha_t N - |\mathbf{a}_t|).    
\end{align}
In this use of the Fenchel Duality Theorem, we note that maximization over policies $\pi$ on the right-hand side of \eqref{eq-8} 
with fixed $\lambda_{1:T}$ can be viewed as as a linear program. More detailed discussion of this standard result can be found in \cite{brown2020index}.

The left-hand side of Equation (\ref{eq-8}) equals $\hat{V}^*(N)$. On the right hand side, for fixed $\lambda_{1: T}$,
\begin{align*}
    \mathbb{E}_{\pi} \sum_{t = 1}^T R_t(\mathbf{s}_t,\mathbf{a}_t) + \lambda_t (\alpha_t N - |\mathbf{a}_t|)
    = \mathbb{E}_{\pi} \sum_{t = 1}^T \sum_{i = 1}^N r_t(s_{t, i}, a_{t, i}) + \lambda_t(\alpha_t - a_{t, i}).
\end{align*}
Since all arms share the same transition kernel, reward function, and distribution over initial state,
\begin{align*}
    \mathbb{E}_{\pi} \sum_{t = 1}^T \sum_{i = 1}^N r_t(s_{t, i}, a_{t, i}) + \lambda_t(\alpha_t - a_{t, i}) = N \ \mathbb{E}_{\pi} \sum_{t = 1}^T r_t(s_{t, 1}, a_{t, 1}) + \lambda_t( \alpha_t - a_{t, 1}).
\end{align*}
So we conclude
\begin{align}\label{eq-9}
  \min_{\lambda_{1: T}} \max_{\pi} \mathbb{E}_{\pi} \sum_{t = 1}^T R_t(\mathbf{s}_t,\mathbf{a}_t) + \lambda_t (\alpha_t N - |\mathbf{a}_t|)
    =  N \min_{\lambda_{1: T}} \max_{\pi} \mathbb{E}_{\pi} \sum_{t = 1}^T r_t(s_{t, 1}, a_{t, 1}) + \lambda_t(\alpha_t - a_{t, 1}).
\end{align}

By using Fenchel Duality again on the one-arm problem,
\begin{align}\label{eq-10}
    \min_{\lambda_{1: T}} \max_{\pi} \mathbb{E}_{\pi} \sum_{t = 1}^T r_t(s_{t, 1}, a_{t, 1}) + \lambda_t(\alpha_t - a_{t, 1}) 
    &= \max_{\pi} \min_{\lambda_{1: T}} \mathbb{E}_{\pi} \sum_{t = 1}^T r_t(s_{t, 1}, a_{t, 1}) + \lambda_t(\alpha_t - a_{t, 1}) \nonumber \\
    &= \hat{V}^*(1).
\end{align}

Summarizing, equations (\ref{eq-8}), (\ref{eq-9}) and (\ref{eq-10}) together imply,
\begin{align*}
    \hat{V}^*(N) = N \ \hat{V}^*(1).
\end{align*}


\subsection{Discussion of the rounding error in budget constraints and initial states}\label{al:2-1}
The original problem \eqref{OP} constrains the number of pulls to $\lfloor \alpha_t N \rfloor$ (almost surely), while the relaxed problem \eqref{KP} constrains this number to $\alpha_t N$ (in expectation).
We think of these differences as ``rounding errors'' in the relaxed problem.
Here we discuss their effect and show that they result in at most a constant difference in the optimal objective value. 

Mathematically speaking, denote
\begin{align*}
    \hat{V}_N^* &= \max_{\pi} \left\{ \mathbb{E}_{\pi} \sum_{t=1}^T R_t(\mathbf{s}_t, \mathbf{a}_t) \Bigg| \mathbb{E}|\mathbf{a}_t| = \alpha_t N, \ \mathbb{E} \sum_{i=1}^N 1(s_{1, i} = s^*) = N, \forall t \in [T] \right\},  \\
    \hat{V}_{N, R}^* &= \max_{\pi} \left\{ \mathbb{E}_{\pi} \sum_{t=1}^T R_t(\mathbf{s}_t, \mathbf{a}_t) \Bigg| \mathbb{E}|\mathbf{a}_t| = \lfloor \alpha_t  N\rfloor, \ \mathbb{E} \sum_{i=1}^N 1(s_{1, i} = s^*) = N, \forall t \in [T] \right\}.
\end{align*}
We claim that $|\hat{V}_{N}^* - \hat{V}_{N, R}^*| \le c$, where $c$ does not depend on $N$. 
The theoretical analysis through the rest of the paper after Lemma~1 compares policy performance against $\hat{V}^*_N$ and shows that this difference is $o(N)$, $O(\sqrt{N})$, or $O(1)$ depending on conditions. 
The fact that $\hat{V}^*_N$ is separated from $\hat{V}_{N,R}$ by at most a constant then implies that the difference in policy performance compared to $\hat{V}^*_{N,R}$ has the same asymptotic dependence on $N$.
This and the fact that $\hat{V}_{N,R}$ is an upper bound on \eqref{OP} 
even when  $\alpha_t N$ are not integers 
provides optimality gaps of $o(N)$, $O(\sqrt{N})$ or $O(1)$ respectively.


The proof of the claim that $|\hat{V}^*_N - \hat{V}^*_{N,R}| \le c$ is 
straightforward. 
As seen from Lemma \ref{fenchel}, there exists a single-arm strategy that pulls $\alpha_t$ arms per period in expectation and achieves objective value $\hat{V}_{1}^*$. Thus, we can pull $N - \max_{t}\ \lceil \frac{1}{\alpha_t} \rceil$ arms according to this strategy and pull each remaining arm with probability $\frac{\lfloor \alpha_t N \rfloor - \alpha_t (N - \max_{t}\ \lceil \frac{1}{\alpha_t} \rceil)}{\max_{t}\ \lceil \frac{1}{\alpha_t} \rceil} \in [0,1]$ at period $t$. Thus, we show
\begin{align*}
    \frac{N - \max_{t}\ \lceil \frac{1}{\alpha_t} \rceil}{N}  \hat{V}_{N}^* - \hat{V}_{N, R}^* \le \ T \max_{s, a, t} r_t(s, a).
\end{align*}
Similarly, we can show
\begin{align*}
    \frac{N - \max_t \ \lceil \frac{1}{\alpha_t} \rceil}{N}  \hat{V}_{N, R}^* - \hat{V}_{N}^* \le \ T \max_{s, a, t} r_t(s, a).
\end{align*}

Combining the above two inequalities 
with the fact that $\hat{V}^*_{N,R}/N$ and $\hat{V}^*_N/N$ are both uniformly bounded by $T \max_{s,a,t} |r_t(s,a)|$
concludes the statement
with $c = 
T(1 + \max_{t}\ \lceil \frac{1}{\alpha_t} \rceil)
\max_{s,a,t} |r_t(s,a)|$.

\subsection{Proof of Lemma \ref{fluid_consistent}}

We prove Lemma \ref{fluid_consistent} by induction on $t$. 
When $t= 1$, that all arms are in state $s^*$ implies
$\frac{Z_1^N}{N} \rightarrow z_1$.
Then, by the definition of fluid consistency,
    $\frac{X_1^N}{N} \rightarrow x_1$.
Thus Lemma \ref{fluid_consistent} holds for $t= 1$.

Now assume Lemma \ref{fluid_consistent} holds for $t$, and we will show it holds for $t + 1$. 
By the definition of fluid consistency, we only need to prove
\begin{align}\label{lemma2_1}
    \frac{Z_{t+1}^N}{N} \rightarrow z_{t+1}.
\end{align}

Recalling our system dynamics,
\begin{equation}\label{dynamic}
    Z^N_{t+1}(s) = \sum_{s' \in S, a \in A} \sum_{i=1}^{N} 1(s_{t, i} = s', a_{t, i} = a, s_{t + 1, i} = s),
\end{equation}
where $1(s_{t, i} = s', a_{t, i} = a, s_{t + 1, i} = s)$ is the indicator function of the event $\{s_{t, i} = s', a_{t, i} = a, s_{t + 1, i} = s\}$,
we only need to show that
\begin{equation}\label{x}
    \frac{1}{N} \sum_{i=1}^{N} 1(s_{t, i} = s', a_{t, i} = a, s_{t + 1, i} = s) \rightarrow x_t(s', a)p(s', a, s).
\end{equation}
since the sum over $s'$ of the right-hand side is $\sum_{s' \in S} x_t(s',a) p(s',a,s) = z_{t+1}(s)$.

If $x_t(s', a) > 0$, then $X^N_t(s', a) \rightarrow \infty$ as $N \rightarrow \infty$ by the induction hypothesis. Thus as $N \rightarrow \infty$, 
\begin{align*}
    \frac{1}{N} \sum_{i=1}^{N} 1(s_{t, i} = s', a_{t, i} = a, s_{t + 1, i} = s) &= \frac{X^N_t(s', a)}{N} \frac{1}{X^N_t(s', a)} \sum_{i=1}^{N}  1(s_{t, i} = s', a_{t, i} = a, s_{t + 1, i} = s) \\
    &\rightarrow x_t(s', a)p(s', a, s),
\end{align*}
by the definition of fluid consistency and the strong law of large numbers.

If $x_t(s', a) = 0$,
\begin{align*}
    \frac{1}{N} \sum_{i=1}^{N} 1(s_{t, i} = s', a_{t, i} = a, s_{t + 1, i} = s) &\leq \frac{X^N_t(s', a)}{N} \rightarrow 0.
\end{align*}

Combining the cases $x_t(s', a) > 0$ and $x_t(s', a) = 0$, equation (\ref{x}) is shown.

To summarize,
\begin{equation*}
    \frac{Z^N_{t+1}(s)}{N}
    = \sum_{s' \in S} \frac{1}{N} \sum_{i=1}^{N} 1(s_{t, i} = s', a_{t, i} = a, s_{t + 1, i} = s)
    \rightarrow \sum_{s' \in S} \sum_{a \in A} x_t(s', a) p(s', a, s)
    = z_{t+1}(s).
\end{equation*}

\subsection{Proof of Lemma \ref{fluid-imply-diffusion}}
This section proves Lemma \ref{fluid-imply-diffusion}. 
When $Z_t^N / N \rightarrow z_t$, we have
\begin{align*}
    \frac{\tilde{Z}_t^N}{\sqrt{N}} = \frac{Z_t^N - N z_t}{N} \rightarrow 0.
\end{align*}

To show $\frac{X_t}{N} \rightarrow x_t$, it is equivalent to show
\begin{align*}
    \frac{\pi_{t, N}(\tilde{Z}_t^N)}{\sqrt{N}} \rightarrow 0.
\end{align*}
Notice by Conditions 1 and 2 in Definition \ref{def-diff},
\begin{align*}
    |\pi_{t, N}(\tilde{Z}_t^N)| \le |\pi_{t, N}(0)| + C_1 |\tilde{Z}_t^N| \le C_2 + C_1 |\tilde{Z}_t^N|.
\end{align*}

Thus,
\begin{align*}
    \frac{\pi_{t, N}(\tilde{Z}_t^N)}{\sqrt{N}} \le \frac{C_2}{\sqrt{N}} + C_1 \frac{|\tilde{Z}_t^N|}{\sqrt{N}} \rightarrow 0.
\end{align*}

\subsection{Proof of Lemma \ref{diffusion_consistent}}
This section proves Lemma \ref{diffusion_consistent}. To begin with, we first state and prove Lemma \ref{le-7}.
\begin{lemma}\label{le-7}
If a policy $\pi$ is diffusion regular and $\tilde{Z}_t^N \rightarrow \tilde{Z}_t^{\infty}$ in distribution, 
then 
$\tilde{X}_t^{N} \rightarrow \tilde{X}_t^{\infty}$
in distribution for some random variable $\tilde{X}_t^{\infty}$.
\end{lemma}

\proof{Proof of Lemma \ref{le-7}}
By the Skorokhod representation Theorem, there exists a probability space $(\Omega, \mathbb{P})$ and a sequence of random variables $\{ \widetilde{Z}_t^N \}_{N}$ and $\widetilde{Z}_{t}^{\infty}$ such that
\begin{align*}
    \widetilde{Z}_t^N &= \tilde{Z}_t^N \text{ and } 
    \widetilde{Z}_t^{\infty} = \tilde{Z}_t^{\infty} \text{ in distribution,}  \\
    \widetilde{Z}_t^N &\rightarrow \widetilde{Z}_t^{\infty} \text{ as $N \rightarrow \infty$}  \quad a.s.
\end{align*}

We will prove convergence in distribution of 
    $\tilde{X}_t^{N} = 
   \tilde{\pi}_{t, N} (\widetilde{Z}_t^N)$
    to 
    $\tilde{X}_t^{\infty} := 
    \tilde{\pi}_{t, \infty} (\widetilde{Z}_t^{\infty})$.
Notice
\begin{align*}
    | \tilde{\pi}_{t, N} (\widetilde{Z}_t^N) -  \tilde{\pi}_{t, \infty} (\widetilde{Z}_t^{\infty})|
    &\leq |\tilde{\pi}_{t, N} (\widetilde{Z}_t^N) -  \tilde{\pi}_{t, N} (\widetilde{Z}_t^{\infty})| +|\tilde{\pi}_{t, N} (\widetilde{Z}_t^{\infty}) -  \tilde{\pi}_{t, \infty} (\widetilde{Z}_t^{\infty})| \\
    &\leq C_1 |\widetilde{Z}_t^N - \widetilde{Z}_t^{\infty}| + |\tilde{\pi}_{t, N} (\widetilde{Z}_t^{\infty}) -  \tilde{\pi}_{t, \infty} (\widetilde{Z}_t^{\infty})|,
\end{align*}
which converges to 0 as $N \rightarrow \infty$ by almost sure convergence of $\widetilde{Z}^N_t$ to $\widetilde{Z}_t^\infty$ and the convergence of $\tilde{\pi}_{t, N}$ to $\tilde{\pi}_{t, \infty}$  required by the fact that $\pi$ is diffusion regular.
\hfill $\Box$


\proof{Proof of Lemma \ref{diffusion_consistent}}
We prove Lemma \ref{diffusion_consistent} by induction on $t$. When $t = 1$, 
    $\tilde{Z}_1^N = 0$ implying $\tilde{Z}_1^\infty=0$.
Then, according to Lemma \ref{le-7}, we know there exists a constant vector $\tilde{X}_1^{\infty}$ s.t.
    $\tilde{X}_1^N \rightarrow \tilde{X}_1^{\infty}$.
Thus, Lemma \ref{diffusion_consistent} holds true for $t = 1$. 

Now assume Lemma \ref{diffusion_consistent} holds for $t$ and we will prove it holds for $ t + 1$. 
It is sufficient to prove there exists a sub-Gaussian random vector $\tilde{Z}_{t+1}^{\infty}$ s.t.
    $\tilde{Z}_{t+1}^{N} \rightarrow \tilde{Z}_{t+1}^{\infty}$
in distribution. 
This is because (1) existence of the limit $\tilde{X}_{t + 1}^{\infty}$ follows from Lemma \ref{le-7} and (2) 
showing $\tilde{Z}^\infty_{t+1}$ is sub-Gaussian implies 
$\tilde{X}_{t+1}^\infty$ is sub-Gaussian because
\begin{equation*}
     |\tilde{X}_{t+1}^{\infty}| 
     = |\tilde{\pi}_{t, \infty}(\tilde{Z}_{t+1}^{\infty})|
     \leq |\tilde{\pi}_{t, \infty}(\tilde{Z}_{t+1}^{\infty}) - \tilde{\pi}_{t, \infty}(0)| + |\tilde{\pi}_{t, \infty}(0)|
     \leq C_1 |\tilde{Z}_{t+1}^{\infty}| + C_2.
\end{equation*}

We prove the existence of $\tilde{Z}_{t+1}^{\infty}$
by constructing an explicit formula for this limit,
\begin{align}\label{intermediate}
\tilde{Z}_{t+1}^{N} \rightarrow \tilde{Z}_{t+1}^{\infty} := M_t + \sum_{s' \in S} \sum_{a \in A} p(s', a, \cdot) \tilde{X}_{t}^{\infty}(s', a)
\end{align}
where $M_t \sim N(0, \Sigma_t)$ is independent of $\tilde{X}_{t}^{\infty}(s', a)$.
The covariance matrix $\Sigma_t$ is defined as
\begin{align*}
    \Sigma_t(s'', s''') = \sum_{s'} \sum_{a \in A} &\ x_t(s', a) \mathrm{Cov}[1(s_{t + 1, 1} = s''), 1(s_{t + 1, 1} = s''') | s_{t, 1} = s', a_{t, 1} = a]
\end{align*}
where $\mathrm{Cov}[1(s_{t + 1, 1} = s''), 1(s_{t + 1, 1} = s''') | s_t = s', a_t = a]$ is the conditional covariance of the indicators of events $\{s_{t + 1, 1} = s''\}$ and $\{s_{t + 1, 1} = s'''\}$ given $s_{t, 1} = s', a_{t, 1} = a$.

Once (\ref{intermediate}) is shown, then the fact that $\tilde{Z}_{t+1}^{\infty}$ is sub-Gaussian follows because $M_t$ and 
$\tilde{X}_{t}^{\infty}$ are both sub-Gaussian.

To prove (\ref{intermediate}), by our system dynamics (\ref{dynamic}) with a vector form,
\begin{equation}
    \label{lemma4-proof-1}
    Z_{t+1}^N = \sum_{s' \in S, a \in A} B_t^N(s', a),
\end{equation}
where the $B_t^N(s', a)$ are conditionally independent (across $s'$ and $a$) multinomial distributions with parameters 
$X_t^N(s', a)$ and $p(s',a) := (p(s', a, s))_{s \in S}$, i.e.,
\begin{equation*}
B_t^N(s', a) \mid X_t^N \sim \mathrm{Multinomial}(X_t^N(s', a), p(s', a)).
\end{equation*}
$B_t^N(s',a)$ is a vector counting the number of arms in each state, among those arms that were previously in state $s'$ and for which we used action $a$.




Recall that $X_t^N(s', a)$ can be decomposed as $N x_t(s', a) + \sqrt{N} \tilde{X}_t^N(s', a)$. According to Lemma \ref{decompose-bin}, there exists two random variables $C_t^N(s', a)$ and $\Delta_t^N(s', a)$, such that
\begin{align}\label{decomp-eq}
    B_t^N(s', a) = C_t^N(s', a) + \Delta_t^N(s', a),
\end{align}
and that, conditionally on $X_t^N(s', a)$, have marginal distributions:
\begin{align*}
    &C_t^N(s', a) \mid X_t^N \sim \mathrm{Multinomial}(N x_t(s', a), p(s', a)),\\
    &\Delta_t^N(s', a)  \mid X_t^N \sim \mathrm{sgn}(\tilde{X}_t^N(s', a)) \mathrm{Multinomial}(\sqrt{N}\ \left| \tilde{X}_t^N(s', a)\right|,\ p(s', a)).
\end{align*}

By \eqref{lemma4-proof-1}, 
\eqref{decomp-eq} and the definition of our diffusion statistic $\tilde{Z}_{t+1}^N$ in terms of $Z_{t+1}^N$,
\begin{align*}
    \tilde{Z}_{t+1}^N = 
    &\frac{1}{\sqrt{N}} \sum_{s' \in S, a \in A}  C_t^N(s', a) - x_t(s', a) p(s', a) N
    + \Delta_t^N(s', a).
\end{align*}

By Lemma \ref{diffusion_tail}, 
\begin{align*}
    \frac{1}{\sqrt N} \Delta_t^N(s', a) - p(s', a)\tilde{X}^N_t(s', a) \rightarrow 0.
\end{align*}

Thus, 
\begin{align*}
    \tilde{Z}_{t+1}^{N}
    &= \left[ \frac{1}{\sqrt{N}} \sum_{s' \in S, a \in A} C_t^N(s', a) - x_t(s', a) p(s', a) N\right] + \left[ \sum_{s' \in S, a \in A}p(s', a)\tilde{X}^N_t(s', a) \right] + \epsilon_N \\
    &= \left[ \frac{1}{\sqrt{N}} \sum_{s' \in S, a \in A} C_t^N(s', a) - x_t(s', a) p(s', a) N\right] + \left[ \sum_{s' \in S, a \in A}p(s', a)\tilde{X}^{\infty}_t(s', a) \right] + \epsilon_N + \epsilon'_N,
\end{align*}
where $\epsilon_N, \epsilon'_N \rightarrow 0$.

The first term satisfies
\begin{align*}
    \frac{1}{\sqrt{N}} \sum_{s' \in S, a \in A} C_t^N(s', a) - x_t(s', a) p(s', a) N \rightarrow N(0, \Sigma_t),
\end{align*}
where $\Sigma_t$ is defined above.
We define $M_t$ to be equal to this limit.
This shows \eqref{intermediate} as claimed.
Although it is not needed for the proof, we observe that because $\tilde{X}_t^\infty$ was constructed to be equal only in distribution to $\lim_N\tilde{X}_t^N$, we are free to construct it so that it is independent of $M_t$.


To summarize, we have shown $\tilde{Z}_{t+1}^{N} \rightarrow \tilde{Z}_{t+1}^{\infty}$ in distribution.
\hfill $\Box$

Here we give the statement and proof of Lemma \ref{decompose-bin} and \ref{diffusion_tail}.

\begin{lemma}\label{decompose-bin}
Let $Y \sim \mathrm{Multinomial}(n,p)$.
Then for a given non-negative integer $m$, there exist random vectors $Y_1$ and $Y_2$ such that
$Y = Y_1 + Y_2$ and
\begin{align*}
    Y_1 \sim \mathrm{Multinomial}(m, p), \ Y_2 \sim \mathrm{sgn}(n - m) \mathrm{Multinomial}(|n - m|, p),
\end{align*}
where $\mathrm{sgn}(\cdot)$ is the sign function.
\end{lemma}

\proof{Proof of Lemma \ref{decompose-bin}}
There exists a sequence of i.i.d random vectors $X_i \sim \mathrm{Multinomial}(1, p)$ s.t. 
$$Y = \sum_{i = 1}^n X_i.$$
If $n > m$, taking $Y_1 = \sum_{i=1}^m X_i, Y_2 = \sum_{j = m + 1}^n X_j$ concludes the proof. If $n \leq m$, taking $Y_1 = \sum_{i=1}^m X_i, Y_2 = -\sum_{j = n + 1}^m X_j$ concludes the proof.
\hfill $\Box$

\begin{lemma}\label{diffusion_tail}
Consider a sequence of random variables $X_1$, $X_2$, ..., $X_N$, ...  converging to $X_{\infty}$ in distribution and a sequence of i.i.d Bernoulli random variable $B_1$, $B_2$, ..., $B_n$, .... with $\mathbb{E} [B_1] = p$ that are also independent of sequence $X_1$, $X_2$, .... Then define
\begin{equation*}
    Y_N = \frac{1}{\sqrt N} \sum_{n=1}^{X_N \sqrt{N}} (B_n - p).
\end{equation*}
Then $Y_N \rightarrow 0$.
\end{lemma}

\proof{Proof of Lemma \ref{diffusion_tail}}
We calculate the characteristic function of $Y_N$ as follows:
\begin{align*}
    \mathbb{E} [\exp(i \lambda Y_N)]
    &= \mathbb{E} [{\mathbb{E} [\exp( i \lambda Y_N) | X_N]]} \\
    &= \mathbb{E} [{\mathbb{E} [\exp( i \lambda \sum_{n=1}^{X_N \sqrt{N}} \frac{1}{\sqrt N}(B_n - p)) | X_N]]} \\
    &= \mathbb{E} [{\mathbb{E} [\exp( i \lambda \frac{B_1 - p}{\sqrt N}) | X_N]^{X_N \sqrt{N}} ]}\\
    &= \mathbb{E} [( p \exp( i \lambda \frac{1-p}{\sqrt N})
    + (1-p) \exp( -i \lambda \frac{p}{\sqrt N}) )^{X_N \sqrt{N}}].
\end{align*}

We have
\begin{align*}
    &(p \exp( i \lambda \frac{1-p}{\sqrt N})
    + (1-p) \exp( -i \lambda \frac{p}{\sqrt N}))^{X_N \sqrt{N}} \\
    &= (p (1 + i \lambda \frac{1-p}{\sqrt N} + O(\frac{1}{N}))
    + (1-p) (1 - i \lambda \frac{p}{\sqrt N} + O(\frac{1}{N})))^{X_N \sqrt{N}} \\
    &= (1 + O(\frac{1}{N}))^{X_N \sqrt{N}}
    \rightarrow 1, \text{ as $N \rightarrow \infty$}.
\end{align*}

We would like to then argue that this almost sure convergence implies convergence of the expectations as well, i.e., that $\mathbb{E}[\exp(i\lambda Y_N)]$ converges to $1$.
To show this we use the dominated convergence theorem and the following bound:
\begin{align*}
|p \exp( i \lambda \frac{1-p}{\sqrt N})
    + (1-p) \exp( -i \lambda \frac{p}{\sqrt N})|
     &\leq p |\exp( i \lambda \frac{1-p}{\sqrt N})|
     + (1-p) |\exp( -i \lambda \frac{p}{\sqrt N})| 
     = p + (1 - p) 
     = 1.
 \end{align*}

Thus $\mathbb{E} [\exp(i \lambda Y_N)] \rightarrow 1$, which implies $Y_N \rightarrow 0$.
\hfill $\Box$

\subsection{Proof of Lemma \ref{moment-bound}}
We only need to prove there exists a constant $C$ s.t. $\mathbb{E}[||\tilde{Z}_t^N||_2^2] \leq C$ for all $t \in [T]$ and $N$. The claim in the lemma for $\tilde{X}_t^N$ follows directly from diffusion regularity. Because
\begin{align*}
    |\tilde{X}_t^N| 
    = |\tilde{\pi}_{t, N}(\tilde{Z}_t^N)|
    \leq |\tilde{\pi}_{t, N}(\tilde{Z}_t^N) - \tilde{\pi}_{t, N}(0)| + |\tilde{\pi}_{t, N}(0)|
    \leq C_1 |\tilde{Z}_t^N| + C_2,
\end{align*}
we have
\begin{align*}
    ||\tilde{X}_t^N||_2^2
    \leq |S| |\tilde{X}_t^N|^2
    \leq |S| |C_1 |\tilde{Z}_t^N| + C_2|^2
    \leq 2|S| C_1^2 |\tilde{Z}_t^N|^2 + 2 |S| C_2^2
    \leq 2|S|^2 C_1^2 ||\tilde{Z}_t^N||_2^2 + 2 |S| C_2^2.
\end{align*}
By taking the expectation,
\begin{align}\label{lemma-5-1}
\mathbb{E}||\tilde{X}_t^N||_2^2 \leq 2|S|^2 C_1^2 \ \mathbb{E}||\tilde{Z}_t^N||_2^2 + 2 |S| C_2^2.    
\end{align}

Similar to the analysis in the proof of Lemma \ref{diffusion_consistent},
\begin{align*}
\tilde{Z}_{t + 1}^N = 
&\frac{1}{\sqrt{N}} \sum_{s' \in S, a \in A} C_t^N(s', a) - x_t(s', a) p(s', a) N + \Delta_t^N(s', a),
\end{align*}
where 
\begin{align*}
    &C_t^N(s', a) \mid X_t^N \sim \mathrm{Multinomial}(N x_t(s', a), p(s', a)),\\
    &\Delta_t^N(s', a) \mid X_t^N \sim \mathrm{sgn}(\tilde{X}_t^N(s', a)) \mathrm{Multinomial}(\sqrt{N} \left| \tilde{X}_t^N(s', a) \right|, p(s', a)).
\end{align*}

Thus, 
\begin{align*}
    ||\tilde{Z}_{t+1}^{N}||_2^2 
    &= \left| \left| \frac{1}{\sqrt{N}} \sum_{s' \in S, a \in A} C_t^N(s', a) - x_t(s', a) p(s', a) N + \Delta_t^N(s', a)\right| \right|_2^2\\
    &\leq 
    2 \left|\left|\frac{1}{\sqrt{N}} \sum_{s' \in S, a \in A} C_t^N(s', a) - x_t(s', a) p(s', a) N\right|\right|_2^2 + 
    2 \left|\left| \frac{1}{\sqrt{N}} \sum_{s' \in S, a \in S} \Delta_t^N(s', a) \right|\right|_2^2.
\end{align*}

Notice from its definition as a multinomial random variable that the absolute value of each component of $\frac{1}{\sqrt{N}} \Delta_t^N(s', a)$ is bounded above by $|\tilde{X}^N_t(s', a)|$. Thus
\begin{align}\label{proof-5-1}
    \left|\left| \frac{1}{\sqrt{N}} \sum_{s' \in S, a \in S} \Delta_t^N(s', a) \right|\right|_2^2 
    \le |S|  \left[ \sum_{s' \in S, a \in S} \left| \tilde{X}_t^N(s', a) \right| \right]^2
    \le 2|S|^2 \sum_{s' \in S, a \in S} \left| \tilde{X}_t^N(s', a) \right|^2 = 2|S|^2 ||\tilde{X}_t^N||_2^2.
\end{align}

On the other hand, noting that 
$C_t^N(s', a) - x_t(s', a) p(s', a) N$ has mean 0 and is independent across different $s,a$ to get the first equality, we have
\begin{align}\label{proof-5-2}
    \mathbb{E}
    \left|\left|\frac{1}{\sqrt{N}} \sum_{s' \in S, a \in A} C_t^N(s', a) - x_t(s', a) p(s', a) N\right|\right|_2^2
    &=\frac{1}{N}
    \sum_{s' \in S, a \in A} \mathbb{E} \left| \left| C_t^N(s', a) - x_t(s', a) p(s', a) N\right| \right|_2^2\\
    &=
    \sum_{s' \in S, a \in A} 
    x_t(s', a) (1 - \sum_{s \in S} p(s', a, s)^2).
\end{align}

Combining inequality (\ref{proof-5-1}) and (\ref{proof-5-2}) together,
\begin{align}\label{lemma-5-2}
    \mathbb{E}||\tilde{Z}_{t+1}^{N}||_2^2 
    \leq 2C_4 + 4|S|^2 \mathbb{E}||\tilde{X}^N_t||_2^2.
\end{align}
where $C_4 := \sum_{s' \in S, a \in A} x_t(s', a) (1 - \sum_{s \in S} p(s', a, s)^2)$.

Recall inequality (\ref{lemma-5-1}) and combine it with (\ref{lemma-5-2}) to obtain,
\begin{align*}
    \mathbb{E}||\tilde{Z}_{t+1}^{N}||_2^2 
    \leq 2C_4 + 8|S|^3 C_2^2 + 8|S|^4 C_2^2C_1^2 \mathbb{E}||\tilde{Z}^N_t||_2^2.
\end{align*}

By $\tilde{Z}^N_1 = 0$ and induction, there exists a constant $C$ s.t. $\mathbb{E}[||\tilde{Z}_t^N||_2^2] \leq C$ for all $t \in [T]$ and $N$.
\hfill $\Box$

\subsection{Proof of Theorem \ref{th:3}}\label{proof-thm-3}
Given a fluid-priority policy $\pi$, we directly check whether the induced map $\tilde{\pi}_{t, N}$ satisfies all three conditions in Definition \ref{def-diff}.

\proof{Verification of Condition 1}
Write the induced map $\tilde{\pi}_{t,N}$ as a collection of maps, 
$(\tilde{\pi}_{t, N}^1, ..., \tilde{\pi}_{t, N}^{|S|})$,
one giving each component.
That is, $\tilde{\pi}_{t,N}(\theta)$ is the vector comprised of 
$(\tilde{\pi}_{t, N}^i(\theta): 1 \le i \le |S|)$.

A direct calculation shows each component function, $\tilde{\pi}_{t, N}^i \ (1 \leq i \leq |S|)$, is continuous, piecewise linear, and has bounded gradients when they exist. Mathematically speaking, there exists a constant $\tilde{C}_1$, s.t., for any $\theta$, any $t$ and any $N$,
\begin{align*}
    | \nabla_\theta \tilde{\pi}_{t, N}^i (\theta) | \leq \tilde{C}_1, \text{\ when $\nabla_\theta \tilde{\pi}_{t, N}^i (\theta)$ exists.}
\end{align*}

For any $\theta_1$ and $\theta_2$, there exists a sequence $(\nu^0, \nu^1, ..., \nu^m)$ lying on the line segment between $\theta_1$ and $\theta_2$, s.t.

1. $\tilde{\pi}_{t, N}^i$ restricted on the line segment between $\nu^j$ and $\nu^{j + 1}$ is linear for $j = 0, 1, ..., m - 1$

2. $\nu^0 = \theta_1$ and $\nu^m = \theta_2$.

Thus
\begin{align*}
    |\tilde{\pi}_{t, N}^i (\theta_1) - \tilde{\pi}_{t, N}^i (\theta_2)| \leq \sum_{j = 0}^{m - 1} |\tilde{\pi}_{t, N}^i (\nu^j) - \tilde{\pi}_{t, N}^i (\nu^{j + 1})|
    \leq \sum_{j = 0}^{m - 1} \tilde{C}_1 |\nu^j - \nu^{j + 1}|
    = \tilde{C}_1 |\theta_1 - \theta_2|.
\end{align*}

So by taking $C_1 = |S| \tilde{C}_1$,
\begin{align*}
    |\tilde{\pi}_{t, N} (\theta_1) - \tilde{\pi}_{t, N} (\theta_2)| \leq \sum_{i = 1}^{|S|} |\tilde{\pi}_{t, N}^i (\theta_1) - \tilde{\pi}_{t, N}^i (\theta_2)|
    \leq \sum_{i=1}^{|S|} \tilde{C}_1 |\theta_1 - \theta_2| = C_1 |\theta_1 - \theta_2|. 
\end{align*}

\proof{Verification of Condition 2}
Direct calculation shows $\tilde{\pi}_{t, N} (0) = 0$.

\proof{Verification of Condition 3}
Direct calculation shows $\tilde{\pi}_{t, \infty} (\tilde{Z}_{t, \infty})$ is a linear mapping. The form of this linear mapping differs across the following three cases. We state the results of detailed calculations here providing these linear forms without including the (tedious) calculations themselves.

Case 1. $C^{0}_t \cup C^{-}_t = \emptyset$: 
\begin{align*}
    \tilde{\pi}_{t, \infty} (\tilde{Z}_{t, \infty})(s, 1) &= \tilde{Z}_{t, \infty}(s), \ \text{for each} \ s \in S.
\end{align*}

Case 2. $C^{0}_t  \neq \emptyset$: 
\begin{align*}
    \tilde{\pi}_{t, \infty} (\tilde{Z}_{t, \infty})(s, 1) &= \tilde{Z}_{t, \infty}(s), \ \text{for each} \ s \in C^{+}_t ;\\
    \tilde{\pi}_{t, \infty} (\tilde{Z}_{t, \infty})(s, 1) &= - \sum_{s' \in C^{+}_t} \tilde{Z}_{t, \infty}(s'), \ \text{for the state $s \in C^{0}_t$ with highest priority-score in $C^{0}_t $};\\
    \tilde{\pi}_{t, \infty} (\tilde{Z}_{t, \infty})(s, 1) &= 0, \ \text{otherwise.}
\end{align*}

Case 3. $C^{0}_t  = \emptyset, \ C^{-}_t  \neq \emptyset$:
\begin{align*}
    \tilde{\pi}_{t, \infty} (\tilde{Z}_{t, \infty})(s, 1) &= \tilde{Z}_{t, \infty}(s), \ \text{for each} \ s \in C^{+}_t;\\
    \tilde{\pi}_{t, \infty} (\tilde{Z}_{t, \infty})(s, 1) &= - \sum_{s' \in C^{+}_t} \tilde{Z}_{t, \infty}(s'), \ \text{for the state $s \in C^{-}_t$ with highest priority-score in $C^{-}_t$};\\
    \tilde{\pi}_{t, \infty} (\tilde{Z}_{t, \infty})(s, 1) &= 0, \ \text{otherwise.}
\end{align*}

To summarize,  we prove the induced map of any fluid-priority policy satisfies all three conditions in Definition \ref{def-diff} and thus any fluid-priority policy is diffusion regular.

\subsection{Proof of Lemma \ref{le-5}}
Direct comparison of Algorithm \ref{al:1} and Algorithm \ref{al:2} justifies Lemma \ref{le-5}.

\subsection{Proof of Lemma \ref{tail}}

Before we prove Lemma \ref{tail}, we prove the following preliminary lemma.
\begin{lemma}\label{one-dimensional-tail}
Suppose the non-degeneracy condition holds. Then there exists constants $\delta > 0$ and $C > 0$ s.t., 
    $\forall \epsilon > 0, t \in [T]$, we have
    \begin{align*}
        \mathbb{P}_{\pi_{R}} \left[ |\tilde{Z}_t^N| \geq \epsilon \sqrt{N} \right] \leq C \exp(- N \delta \epsilon^2).
    \end{align*}
\end{lemma}

\proof{Proof of Lemma \ref{one-dimensional-tail}}
We will show that, for $0 \le t \le T - 1$,
\begin{align}\label{dynmc}
    \mathbb{P}_{\pi_{R}} \left[ |\tilde{Z}_{t + 1}^N| \geq \epsilon \sqrt{N} \right] \le 4|S|^2 \exp(-\frac{\epsilon^2 N}{8|S|^4}) + 2|S|^2 \mathbb{P}_{\pi_{R}}\left[|\tilde{Z}_{t}^N| \geq \frac{\epsilon \sqrt{N}}{4|S|^2}\right].
\end{align}
The above inequality and the observation $\tilde{Z}_{1}^N = 0$ would complete Lemma \ref{one-dimensional-tail}. So in the remainder of this proof, we show inequality (\ref{dynmc}) holds true.

By a union bound and the fact that $|\tilde{Z}_{t+1}^N|\ge \epsilon \sqrt{N}$ implies $|\tilde{Z}_{t+1}^N(s)|\ge \epsilon \sqrt{N} / |S|$ for at least one $s$,
\begin{align*}
    \mathbb{P}_{\pi_{R}} \left[ |\tilde{Z}_{t + 1}^N| \geq \epsilon \sqrt{N} \right] \le \sum_{s \in S} \mathbb{P}_{\pi_{R}} \left[ |\tilde{Z}_{t + 1}^N(s)| \geq \frac{\epsilon \sqrt{N}}{|S|} \right].
\end{align*}
Thus, we only need to show, for any $s \in S$,
\begin{align}\label{one-state-analysis}
    \mathbb{P}_{\pi_{R}} \left[ |\tilde{Z}_{t + 1}^N(s)| \geq \frac{\epsilon \sqrt{N}}{|S|} \right] \le 4|S| \exp(-\frac{\epsilon^2 N}{8|S|^4}) + 2|S| \mathbb{P}_{\pi_{R}}(|\tilde{Z}_{t}^N| \geq \frac{\epsilon \sqrt{N}}{4|S|^2}).
\end{align}

Following a similar approach to the proof of Lemma \ref{diffusion_consistent}, we first write our system dynamics in a vector form:
\begin{equation}\label{lemma12-1}
    Z_{t+1}^N = \sum_{s' \in S, a \in A} B_t^N(s', a),
\end{equation}
where the $B_t^N(s', a)$ are conditionally independent (across $s'$ and $a$) multinomial distributions with parameters 
$X_t^N(s', a)$ and $p(s',a) := (p(s', a, s))_{s \in S}$, i.e.,
\begin{equation*}
B_t^N(s', a) \mid X_t^N \sim \mathrm{Multinomial}(X_t^N(s', a), p(s', a)).
\end{equation*}
$B_t^N(s',a)$ is a vector counting the number of arms in each state, among those arms that were previously in state $s'$ and for which we used action $a$. We use $B_t^N(s',a, s)$ to denote component $s$ of $B_t^N(s',a)$, i.e. the number of arms that were previously in state $s'$, for which we used action $a$, and which transitioned to state $s$.

Recall the definition of our diffusion statistic $\tilde{Z}_{t+1}^N$ and combine it with equation (\ref{lemma12-1}),
\begin{align*}
    \tilde{Z}_{t+1}^{N}(s) 
    &=  \left[ \sum_{s' \in S, a \in A} \frac{1}{\sqrt{N}} (B_t^N(s', a, s) - p(s', a, s) X_t^N(s', a)) \right] +
    \left[ \sum_{s' \in S, a \in A} p(s', a, s) \tilde{X}_t^N(s', a) \right].
\end{align*}

So we have
\begin{align*}
    \mathbb{P}_{\pi_{R}} \left[ |\tilde{Z}_{t + 1}^N(s)| \geq \frac{\epsilon \sqrt{N}}{|S|} \right] 
    &\le \mathbb{P}_{\pi_{R}} \left[ \left|\sum_{s' \in S, a \in A} \frac{1}{\sqrt{N}} (B_t^N(s', a, s) - p(s', a, s) X_t^N(s', a))\right| \geq \frac{\epsilon \sqrt{N}}{2|S|} \right]\\
    &+ \mathbb{P}_{\pi_{R}} \left[ \left|\sum_{s' \in S, a \in A} p(s', a, s) \tilde{X}_t^N(s', a)\right| \geq \frac{\epsilon \sqrt{N}}{2|S|} \right].
\end{align*}

Notice
\begin{align*}
    &\mathbb{P}_{\pi_{R}} \left[ \left|\sum_{s' \in S, a \in A} \frac{1}{\sqrt{N}} (B_t^N(s', a, s) - p(s', a, s) X_t^N(s', a))\right| \geq \frac{\epsilon \sqrt{N}}{2|S|} \right]\\
    &\leq \sum_{s' \in S, a \in A} \mathbb{P}_{\pi_{R}}\left[ 
    \frac{1}{\sqrt{N}} \left| B_t^N(s', a, s) - p(s', a, s) X_t^N(s', a) \right|
    \geq \frac{\epsilon \sqrt{N}}{4 |S|^2}  \right].
\end{align*}
By Hoeffding's inequality, we have
\begin{align*}
    &\mathbb{P}_{\pi_{R}}\left[ \frac{1}{\sqrt{N}} \left| B_t^N(s', a, s) - p(s', a, s) X_t^N(s', a) \right| \geq \frac{\epsilon \sqrt{N}}{4 |S|^2} \right] \leq 2 \exp\left(-\frac{\epsilon^2 N^2}{8 |S|^4 |X_t^N(s', a)|}\right) \leq 2 \exp\left(-\frac{\epsilon^2 N}{8 |S|^4}\right).
\end{align*}
Combining the above inequalities together,
\begin{align}\label{fluid-pull}
\mathbb{P}_{\pi_{R}} \left[ \left|\sum_{s' \in S, a \in A} \frac{1}{\sqrt{N}} (B_t^N(s', a, s) - p(s', a, s) X_t^N(s', a))\right| \geq \frac{\epsilon \sqrt{N}}{2|S|} \right]
\leq 4|S| \exp(-\frac{\epsilon^2 N}{8 |S|^4}).
\end{align}

On the other hand, combining the bound
\begin{align*}
    \left|\sum_{s' \in S, a \in A} p(s', a, s) \tilde{X}_t^N(s', a)\right| \le \sum_{s' \in S, a \in A} |\tilde{X}_t^N(s', a)|
\end{align*}
with a union bound, we have
\begin{align*}
    \mathbb{P}_{\pi_{R}} \left[ \left|\sum_{s' \in S, a \in A} p(s', a, s) \tilde{X}_t^N(s', a)\right| \geq \frac{\epsilon \sqrt{N}}{2|S|} \right] 
    &\le \sum_{s' \in S, a \in A} \mathbb{P}_{\pi_{R}} \left[  \left| \tilde{X}_t^N(s', a)\right| \geq \frac{\epsilon \sqrt{N}}{4|S|^2} \right].
\end{align*}
Analysis similar to Condition 3 in \S\ref{proof-thm-3} shows that $\left| \tilde{X}_t^N(s', a)\right| \le |\tilde{Z}_t^N|$ for any $s \in S$. Thus,
\begin{align}\label{diffusion-pull}
    \mathbb{P}_{\pi_{R}} \left[ \left|\sum_{s' \in S, a \in A} p(s', a, s) \tilde{X}_t^N(s', a)\right| \geq \frac{\epsilon \sqrt{N}}{2|S|} \right] 
    &\le 2|S| \mathbb{P}_{\pi_{R}} \left[  |\tilde{Z}_t^N| \geq \frac{\epsilon \sqrt{N}}{4|S|^2} \right].
\end{align}

Combining inequality (\ref{fluid-pull}) and (\ref{diffusion-pull}) together implies inequality (\ref{one-state-analysis}), concluding the proof.
\hfill $\Box$

Now we can prove Lemma \ref{tail}. 
\proof{Proof of Lemma \ref{tail}}
Let $\Omega_t := \Delta_1 \cap \Delta_2 \cap\ldots\cap\Delta_t$ and let $\Omega_t^c$ denote its complement.
First we notice,
\begin{align*}
    \mathbb{P}_{\pi_{F}}(\Delta_{t+1}^c)
    &= \mathbb{P}_{\pi_{F}}(\Omega_t \cap \Delta_{t+1}^c) + \mathbb{P}_{\pi_{F}}(\Omega_t^c \cap \Delta_{t+1}^c) \\
    &= \mathbb{P}_{\pi_{R}}(\Omega_t \cap \Delta^c_{t+1}) + \mathbb{P}_{\pi_{F}}(\Omega_t^c \cap \Delta_{t+1}^c) \\
    &\leq \mathbb{P}_{\pi_{R}} (\Delta^c_{t+1}) +  \mathbb{P}_{\pi_{F}}(\Omega_t^c) \\
    &\leq \mathbb{P}_{\pi_{R}} (\Delta^c_{t+1}) +  \sum_{k = 1}^t \mathbb{P}_{\pi_{F}}(\Delta_k^c).
\end{align*}

We will use this recursive expression show that $\mathbb{P}_{\pi_F}(\Delta_t^c) \le L \exp(-\delta N)$ by induction on $t$.
The base case, $t=1$, follows immediately from $\mathbb{P}_{\pi_{R}}(\Delta_1^c) = 0$. Thus, it is sufficient to prove there exists constants $\delta > 0$ and $L$, s.t. for all $t$,
\begin{align}\label{tail_pi}
    \mathbb{P}_{\pi_{R}}(\Delta_t^c) \leq L \exp(-\delta N).
\end{align}

We rewrite $\Delta_t^c$ in terms of $\tilde{Z}_t^N$,
by first noting that there are two ways to have a budget violation event $\Delta_t^c$. The first arises when the number of arms available to pull in fluid-active and fluid-neutral states,
$\sum_{s \in C_t^0 \cup C_t^+} Z_t^N(s)$, falls below the number of arms that the optimal occupation measure plans to pull $N\sum_{s\in C_t^0 \cup C_t^+} x_t(s,1)$, where we note that the optimal occupation measure never pulls arms in $C_t^-$. We define our diffusion statistics $\tilde{Z}_t^N$ by subtracting $(x_t(s,0) + x_t(s,1))N$ from $Z_t^N$ and dividing the difference by $\sqrt{N}$, and so the following conditions are all equivalent:
\begin{align*}
N\sum_{s \in C_t^0 \cup C_t^+} x_t(s,1) &> \sum_{s \in C_t^0 \cup C_t^+} Z_t^N(s)  \\
-N\sum_{s \in C_t^0 \cup C_t^+} x_t(s,0) &> \sum_{s \in C_t^0 \cup C_t^+} Z_t^N(s) - N(x_t(s,0) + x_t(s,1)),\\
-\sqrt{N}\sum_{s \in C_t^0 \cup C_t^+} x_t(s,0) &> \sum_{s \in C_t^0 \cup C_t^+} \tilde{Z}_t^N(s).
\end{align*}
Moreover, optimal occupation measures set $x_t(s,0)=0$ for $s\in C_t^+$. Thus, the conditions above are equivalent to
\begin{equation*}
-\sqrt{N}\sum_{s \in C_t^0} x_t(s,0) > \sum_{s \in C_t^0 \cup C_t^+} \tilde{Z}_t^N(s).
\end{equation*}

The other way in which we can have a budget violation is to have the number of arms available to idle in fluid-inactive and fluid-neutral states fall below the number of arms that the optimal occupation measure plans to idle, $N\sum_{s \in C_t^0 \cup C_t^-} x_t(s,0)$. By a similar sequence of computations, this occurs if and only if
\begin{equation*}
\sum_{s \in C_t^+} \tilde{Z}^N_t(s) > \sqrt{N} \sum_{s\in C_t^0} x_t(s,1)
\end{equation*}

Thus,
\begin{align*}
    \Delta_t^c 
    &= \Big\{ -\sqrt{N} \sum_{s \in C^{0}_t} x_t(s, 0) > \sum_{s \in C^{0}_t \cup C^+_t} \tilde{Z}_t^N(s) \Big\} \bigcup \Big\{ \sum_{s \in C^{+}_t} \tilde{Z}_t^N(s) > \sqrt{N} \sum_{s \in C^{0}_t} x_t(s, 1) \Big\}.
\end{align*}

Thus we have,
\begin{align*}
    \mathbb{P}_{\pi_{R}}(\Delta_t^c) 
    &\leq \mathbb{P}_{\pi_{R}}\left[-\sqrt{N} \sum_{s \in C^{0}_t} x_t(s, 0) > \sum_{s \in C^{0}_t \cup C^+_t} \tilde{Z}_t^N(s_t)
    \right] 
    + \mathbb{P}_{\pi_{R}}\left[\sum_{s \in C^{+}_t} \tilde{Z}_t^N(s) > \sqrt{N} \sum_{s \in C^{0}_t} x_t(s, 1)\right]\\
    &\leq \sum_{s \in C^{0}_t \bigcup C^{+}_t} \mathbb{P}_{\pi_{R}}\left[|\tilde{Z}_t^N(s)| > \sqrt{N} \frac{\sum_{s \in C^{0}_t} x_t(s, 0)}{|S|}\right] + \sum_{s \in C^{0}_t} \mathbb{P}_{\pi_{R}}\left[|\tilde{Z}_t^N(s)| > \sqrt{N} \frac{\sum_{s \in C^{0}_t} x_t(s, 1)}{|S|}\right]
\end{align*}
Using Lemma \ref{one-dimensional-tail}, it is easy to see inequality (\ref{tail_pi}) holds.

\subsection{Proof of Lemma \ref{le-6}}
By the Fenchel Duality Theorem \citep{rockafellar1970convex}, there exists $\lambda^*_{1:T} = (\lambda_1^*, \lambda_2*, ..., \lambda_T^*)$ s.t.
\begin{align}\label{eq-16}
    \hat{V}_1^* = \max_{\pi} \mathbb{E}_{\pi} \sum_{t = 1}^T r_t(s_{t, 1}, a_{t, 1}) + \lambda_t^*( \alpha_t - a_{t, 1})
\end{align}
where here the maximum is taken over all policies, not just those satisfying the budget constraint $\mathbb{E}_\pi |a_{t,1}| = \alpha_t$.

Following a dynamic programming argument, we define the value function $V_t$ on states and the Q-factor $Q_t$ on state-action pairs recursively as
\begin{align*}
    Q_t(s, a) &= r_t(s, a) - \lambda_t^* a + \sum_{s' \in S} p_t(s, 1, s') V_{t + 1}(s'),\\
    V_t(s) &= \max \{ Q_t(s, 0), Q_t(s, 1)\}.
\end{align*}
for $0 \leq t \leq T$ with $V_{T + 1}(s) = 0$ for all $s \in S$. Thus, we can classify states into three disjoint sets
\begin{align*}
    \text{Must-Pull}_t &= \{s \in S | Q_t(s, 0) < Q\},\\
    \text{Indifferent}_t &= \{s \in S | Q_t(s, 0) = Q_t(s, 1) \},\\
    \text{Never-Pull}_t &= \{s \in S | Q_t(s, 0) > Q_t(s, 1) \}.
\end{align*}

A policy is optimal for \eqref{eq-16} if and only if it satisfies these two conditions for each $t$:
\begin{itemize}
\item It pulls all arms whose states are in $\text{Must-Pull}_t$.
\item It never pulls any arms whose states are in $\text{Never-Pull}_t$.
\end{itemize}
It can behave arbitrarily for arms whose states are in $\text{Indifferent}_t$.


Any optimal occupation measure $(x_t(s, a))_{t \in [T], s \in S, a \in A}$ achieves $\hat{V}_1^*$ and so must correspond to an optimal policy.
Thus
\begin{align*}
    \text{Must-Pull}_t \subseteq C^{+}_t, \ \text{Never-Pull}_t \subseteq C^{-}_t.
\end{align*}

Any budget-relaxed fluid-priority policy $\pi_{R}$ 
pulls an arm whenever its state is in $C^{+}_t$ and lets an arm idle whenever its state is in $C^{-}_t$. 
Thus, it is optimal for \eqref{eq-16} and
\begin{align*}
    V_N(\pi_R) + \sum_{t=1}^T \lambda_t^* (\alpha_t N - \mathbb{E}_{\pi_R}[|\textbf{a}_t|]) = \hat{V}_N^*.
\end{align*}
Using the fact that 
$|\alpha_t N - \mathbb{E}_{\pi_R}[\textbf{a}_t]|$ is bounded above by the probability of a budget violation event times a bound $N$ on the maximum size of a budget violation, as well as Lemma \ref{tail},
\begin{align*}
    |\alpha_t N - \mathbb{E}_{\pi_R}[|\textbf{a}_t|]| 
    \le N \mathbb{P}_{\pi_R}(\Delta_t^c)
    \le N L \exp(-\delta N)
    \le  m, \ \forall t,
\end{align*}
where $m$ is a constant not depending on $N$.

Thus,
\begin{align*}
    |V_N(\pi_R) - \hat{V}_N^*| 
    &= |\sum_{t=1}^T \lambda_t^* (\alpha_t N - \mathbb{E}_{\pi_R}[|\textbf{a}_t|])| \le \sum_{t=1}^T |\lambda_t^*| |\alpha_t N - \mathbb{E}_{\pi_R}[|\textbf{a}_t|])| \le m \sum_{t=1}^T |\lambda_t^*|.
\end{align*}

\subsection{Proof for Proposition \ref{prop-1}}
For any index policy $\pi_I$, by the strong law of large numbers,
$X_t^N(s,a)/N$ converges as $N\to\infty$ to a quantity that we denote $x_{I,t}(s,a)$ and refer to as the occupation measure of the index policy.


We argue that $x_{I,t}(s,a)$ has at most one state $s$ for each $t$ satisfying both $x_{I,t}(s,0)>0$ and $x_{I,t}(s,1)>0$.
To see this, first recall that index policies use a strict priority order over states, pulling all arms in states higher in the priority order before pulling any arms in lower states.
Then define for each state $s$ and time $t$ the following quantities:
\begin{itemize}
\item Let $P(s)$ denote the set of states that have equal or higher priority to $s$ according to the index policy.
\item Let $L_{t}^N(s)$ denote the number of arms whose states have equal or higher priority than $s$ and that are not pulled.
\item Let $M_t^N(s)$ denote the number of arms pulled whose states have priority strictly lower than $s$.
\end{itemize}
By the mechanics of an index policy's decisions, we either have $L_t^N(s)= 0$, i.e., we pull all of the arms whose states have equal or higher priority to $s$, or $M_t^N(s) = 0$, i.e., we pull no arms whose states have priority strictly lower than $s$ 

Then, taking the limit as $N\to\infty$ and using the strong law of large numbers, we have
\begin{align*}
 0 = \lim_{N\to\infty} \frac{L_t^N(s)M_t^N(s)}{N^2}
 = \left(\sum_{s' \in P(s)} x_{I,t}(s',0)\right)\left(\sum_{s' \notin P(s)} x_{I,t}(s',1)\right),
\end{align*}

This then implies that there is a unique $s$ such that 
$x_{I,t}(s',0)=1$ for $s' \in P(s) \setminus \{s\}$
and 
$x_{I,t}(s',1)=0$ for all $s'\notin P(s)$.
That is, states that have strictly higher priority than $s$ are always pulled in the fluid limit, while states that have strictly lower priority than $s$ are never pulled in this limit. 

Now, since any index policy meeting the condition of the proposition has $\hat{V}_N^* - V_{N}(\pi_{I})$ bounded above by a constant, this index policy's occupation measure $x_{I,t}$ solves Problem \ref{LP}.
We then construct a fluid-priority policy to match this index policy.

First, we note that the set of fluid-active states for the optimal occupation measure $x_{I,t}$ are those with $x_{I,t}(s,0)=0$ and that the index policy ranks these above all other states. We take the priority score used by our fluid priority policy to rank these fluid-active states among themselves in the same way as the index policy.

Second, the set of fluid-inactive states for $x_{I,t}$ are those with $x_{I,t}(s,1)=0$. The index policy ranks these below all other states. Again, we take the priority score used by our fluid priority policy to rank these fluid-inactive states in the same way as the index policy.

Third, the at most one state with $x_{I,t}(s,0)>0$ and $x_{I,t}(s,1)>0$ is a fluid-neutral state, and it is ranked by the index policy below the fluid-active states and above the fluid-inactive states.

Because our fluid-priority policy's priority score matches the index policy's prioritizations on fluid-active and fluid-inactive states, and its prioritizations also across categories (fluid-active, fluid-neutral, fluid-inactive) match those of the index policy, our fluid-priority policy is the same as the index policy.



\subsection{Proof for Proposition \ref{ucb_ts_regret}}
The proof is similar for both UCB and Thompson Sampling policy. We only show the proof for UCB here.

Under the UCB policy, there exists $z^{UCB}_t(s)$ and $x^{UCB}_t(s, a)$ which are feasible for the LP (\ref{LP}) and satisfy
\begin{align*}
    \frac{Z_t^N(s)}{N} \rightarrow z_t^{UCB}(s), \frac{X_t^N(s, a)}{N} \rightarrow x_t^{UCB}(s, a).
\end{align*}

So we have
\begin{align*}
    \frac{R_{\pi^{UCB}}(N)}{N} \rightarrow &\sum_{t = 1}^T \sum_{s \in S} \sum_{a \in A} r_t(s, a) x^{UCB}_t(s, a).
\end{align*}

UCB policy is an index policy. Thus, the occupation measure $x_t^{UCB}(s, a)$ can be calculated via forward propagation. Numerically, we can verify $x_t^{UCB}(s, a)$ is not an optimal solution for LP (\ref{LP}) under $T=15$ and $T=20$.



\subsection{Discussion of policies in previous literature}\label{other-policies}
In this section, we show the power of 
the techniques developed in \S\ref{oN} and \ref{sec-4} by applying them to policies proposed by previous literature to demonstrate theoretical guarantees from that literature can be seen as consequences of our results.
Specifically, we observe that the Randomized Assignment Control (RAC) policy proposed by \cite{zayas2019asymptotically} is fluid-consistent, thus achieving an $o(N)$ optimality gap.
The policy proposed by \cite{hu2017asymptotically} and the ``optimal Lagrangian index policy'' proposed by \cite{brown2020index} are diffusion-regular, thus achieving $O(\sqrt{N})$ optimality gaps.

\subsubsection{\cite{zayas2019asymptotically} achieves $o(N)$ optimality gap}
This section shows the RAC policy proposed by \cite{zayas2019asymptotically} achieves an $o(N)$ optimality gap. To start with, let us first describe the RAC policy. Although \cite{zayas2019asymptotically} defines
RAC policy in settings more general than the binary-action bandit (referring to their more general problem setting as a ``multi-action bandit''), we only focus on the binary bandit here.

Similar to our approach, \cite{zayas2019asymptotically} first solves the linear programming relaxation (\ref{P}) and then fetch an optimal occupation measure $\{x_t(s, a)\}_{t \in [T], s \in S, a \in A}$. Then, based on the occupation measure, an activation probability is defined for each state $s$ at period $t$:

\[
  q_t(s) =
  \begin{cases}
  \frac{x_t(s, 1)}{z_t(s)}, & \text{if $z_t(s) > 0$}; \\
   \\
  0, & \text{if $z_t(s) = 0$}.
  \end{cases}
\]

Then when deciding which arm to pull at period $t$ under the RAC policy, we first randomly choose an arm that has not been chosen in this period. If the arm's state is $s$, then we randomly generate a Bernoulli variable with mean $q_t(s)$. If this random realization is $0$ or there is no remaining budget, idle the arm; otherwise, activate the arm. Repeat this process until no budget remains in the period.

Direct computation and the strong law of large numbers show that the RAC policy is fluid consistent. Thus, it achieves an $o(N)$ optimality gap.

\subsubsection{\cite{hu2017asymptotically} and \cite{brown2020index} achieve $O(\sqrt{N})$ optimality gaps}
The policies proposed by \cite{hu2017asymptotically} and \cite{brown2020index} are very similar. 
Thus we only discuss \cite{brown2020index}'s policy here. The analysis for \cite{hu2017asymptotically}'s policy can be generalized without any essential difficulty.

To start with, we first describe the ``optimal Lagrangian index policy'' proposed by \cite{brown2020index}. Similar to our approach, \cite{brown2020index} first solves the linear programming relaxation (\ref{P}) and fetches an optimal occupation measure $\{x_{t}(s, a)\}_{t \in [T], s \in S, a \in A}$, which is used to do ``tie-breaking'' discussed later. While solving the relaxed problem using the Simplex method \citep{nash2000dantzig}, as a byproduct, its dual problem
\begin{align*}
    \min_{\lambda_{1: T}} \max_{\pi} \mathbb{E}_{\pi} \sum_{t = 1}^T r_t(s_t, a_t) + \lambda_t (\alpha_t N - a_t).
\end{align*}
is also solved, which yields optimal Lagrange multipliers $\{ \lambda_t^* \}_{t = 1}^T$. 

Then, following a dynamic programming argument, the value function $V_t$ on states and the Q-factor $Q_t$ on state-action pairs are defined as
\begin{align*}
    Q_t(s, a) &= r_t(s, a) - \lambda_t^*a + \sum_{s' \in S} p_t(s, 1, s') V_{t + 1}(s'),\\
    V_t(s) &= \max \{Q_t(s, 0), Q_t(s, 1)\}
\end{align*}
for $0 \le t \le T$ with $V_{T + 1}(s) = 0$ for all $s \in S$. Finally the index of a state $s$ at period $t$ is defined as 
\begin{align*}
\textnormal{Index(s) }= Q_t(s, 1) - Q_t(s, 0).
\end{align*}

When deciding which arm to pull, arms are activated from high index to low index until no budget remains. When there is a tie, i.e., some states share the same index value, the number of arms activated from a state is proportional to its occupation measure. More details can be found in \cite{brown2020index} Section 4.

Now we show the optimal Lagrangian index policy is diffusion regular. First of all, we can show its associated map $\hat{\pi}_{t, N}$ is a piece-wise linear map, thus satisfying Condition 1 in Definition \ref{def-diff}. Second, we can show $\hat{\pi}_{t, N}(0) = 0$, thus satisfying Condition 2. As a piece-wise linear map, we can also show that $\hat{\pi}_{t, N}$ converges as $N \rightarrow \infty$. Thus, Condition 3 is satisfied.
